\documentclass[12pt]{article}
\title{One-cocycle invariants for closed braids}
\author{Thomas Fiedler}
\usepackage{amssymb, amsfonts, epsfig}
\begin{document}
\newtheorem{proposition}{Proposition}
\newtheorem{lemma}{Lemma}
\newtheorem{theorem}{Theorem}
\newtheorem{definition}{Definition}
\newtheorem{remark}{Remark}
\newtheorem{example}{Example}
\newtheorem{corollary}{Corollary}
\newtheorem{observation}{Observation}
\maketitle
\begin{center}
{\em to S\'everine}
\end{center}
\begin{abstract}

We introduce new polynomial  isotopy invariants for closed braids. They are constructed as polynomial valued {\em Gauss diagram 1-cocycles} evaluated on the full rotation of the closed braid $\hat \beta$ around the core of the corresponding solid torus. They can be calculated with polynomial complexity with respect to the braid length and their derivatives evaluated at $x=1$ are finite type invariants of closed braids.

Let the solid torus V be standardly embedded in the 3-sphere and let L be the core of the complementary solid torus 
$S^3\setminus V$. We give examples which show that a natural refinement of our invariants can detect (even with linear complexity with respect to the braid length if the number of strands is fixed, and with quadratic complexity if it is not fixed) the non-invertibility of  the 2-component link $\hat \beta \cup L\hookrightarrow S^3$, what quantum invariants fail to do. 

\footnote{2000 {\em Mathematics Subject Classification\/}: 57M25.{\em Key words and phrases\/}: Closed braids, one-cocycle polynomials, Vassiliev invariants, character invariants, non-invertibility of links}
\end{abstract}
\tableofcontents

\section{Introduction}

This paper is an improved and shortened version of the preprint \cite{F06}. Moreover, it contains a lift of the integer valued 1-cocycles of \cite{F06} to polynomial valued 1-cocycles. They can be calculated with the {\em same complexity} as before but they are no longer of finite type.

As well known, the isotopy problem for 
closed braids in the solid torus reduces to the conjugacy problem in braid groups (see e.g. \cite{M78}).
 The latter problem is solved, but in general only with exponential complexity with respect to the braid length 
(see \cite{G} and \cite{B-B04}).
It is therefore interesting to construct invariants which distinguish conjugacy classes of braids and which are 
calculable with polynomial complexity. Finite type invariants for knots in the solid torus are an example of such invariants 
(see \cite{V90}, \cite{V01}, \cite{BN}, \cite{G96}, \cite{F01} and references therein). We introduce a new approach in this paper, which is based on the fact that each component of the topological moduli space of closed braids in the solid torus has a non-trivial first homology group.

{\em The main result of this paper (Theorems 4,5,6 and 7) is the construction of two new classes of invariants for closed braids, which are calculable in polynomial time:  {\em one-cocycle polynomials and character invariants}. Moreover, our new invariants are more related to geometric invariants of braids (as entropy and simplicial volume of its mapping torus) than are the usual invariants (e.g. the HOMFLYPT polynomial for links in the solid torus, compare \cite {Tu} and \cite{G96} for its definition): if for an irreducible braid one of our invariants is non-trivial then the braid is pseudo-Anosov, and hence both its entropy and its hyperbolic volume are non-trivial.}\vspace{0.4cm}

We give now a brief outline of our approach:  let $\hat \beta\subset V$ be a closed braid (i.e. it intersects each disc in the fibration $V=S^1 \times D^2$ transversely) and such that $\hat \beta$ is a knot. We fix a projection 
$pr : V \to S^1\times I$. Let $M_n$ be the  infinite dimensional space of all closed n-braids in $V$ and which are knots. The space $M_n$ has a natural stratification with respect to $pr$.  Let $\sum$ be the discriminant in $M_n$ which consists of all non-generic diagrams of closed braids. The components of $M_n \setminus \sum$ consist of the ordinary diagrams of closed n-braids. The strata $\sum^{(1)}$ of codimension 1 are just the braid diagrams which have in the projection $pr$ in addition either exactly one ordinary triple point or exactly one ordinary autotangency. We call the corresponding strata $\sum^{(1)}_{tri}$ and respectively $\sum^{(1)}_{tan}$.
First, we associate to a closed braid in a canonical way a loop in $M_n$, called $rot(\hat \beta)$, namely the loop induced by the rotation of the solid torus around its core. This loop can be seen as a generalization of {\em Gramain's loop} for long knots, compare \cite{H}, \cite{F-K}. It has a very nice combinatorial description, namely it is just the pushing of a full-twist $\Delta^2$ through the closed braid. We associate to the canonical loop an oriented singular link in a thickened torus by tracking the crossings of the closed braids in the loop. This link is called the {\em trace graph \/}
and it is denoted by $TL(rot(\hat \beta))$. All singularities of $TL(rot(\hat \beta))$ are ordinary triple points. These triple points 
correspond exactly to the intersections of $rot(\hat \beta)$ with $\sum^{(1)}_{tri}$.
There is a natural coorientation on $\sum^{(1)}_{tri}$ and, hence, each triple point in $TL(rot(\hat \beta))$ has a sign.
To each triple point corresponds a diagram of a closed braid which has just ordinary crossings and exactly one triple crossing.
We use the position of the ordinary crossings with respect to the triple crossing in the {\em Gauss diagram \/}
in order to construct monomials as {\em weights \/} for the triple crossings. 

{\em We associate then to each generic loop $\gamma \subset M_n$  polynomial invariants $\Gamma(\gamma) \in \mathbb{Z} [x,x^{-1}]$ which depend only on the homology class of $\gamma \subset M_n$, called {\em one-cocycle polynomials}, by summing over the weighted intersection numbers of $\gamma$ with $\sum^{(1)}_{tri}$.}

 We show that in particular $\Gamma(\hat \beta) = \Gamma( rot(\hat \beta))$ is not of finite type but that 
$d(\Gamma(\hat \beta))/dx$ evaluated at  $x=1$ is a integer valued finite type invariant of $\hat \beta$. We call them {\em finite type one-cocycle invariants}. As all finite type invariants, the finite type one-cocycle invariants have a natural degree. However, our invariants induce a new filtration on the space of all finite type invariants for closed braids, as shows the following fact:
Let $\hat \beta$ be a closed n-braid which is a knot and let c be the word length of $\beta \in B_n$ (with respect to the standard generators of $B_n$, or in other words, the number of crossings of the diagram with respect to $pr$).
Then all finite type one-cocycle invariants of degree d vanish for $\hat \beta$  if   $d \geq  c + n^2 - n - 1$.

One has to compare this with the following well known fact: already the trefoil has non-trivial finite type invariants 
of arbitrary high degree. Consequently, the finite type one-cocycle invariants define a new natural filtered subspace in the filtered space of all finite type invariants for those knots in the solid torus, which are closed braids.

Our method is based on Thom-Mather singularity theory together with combinatorics. We give extremely simple new solutions of the tetrahedron equation in the case of closed braids, without solving a big system of equations. (They satisfies  automatically the {\em marked 4T-relations}, compare \cite{G96} and also \cite{F01}.)\vspace{0.4cm}

However, it turns out that the one-cocycle polynomials can be refined considerably.
When we deform $\hat \beta$ in V by a generic isotopy then $rot(\hat \beta)$ in $M_n$ deforms by a generic 
homotopy.
The following is our key observation for the refinement (see Section 2.4).
\begin{observation}
The loop $rot(\hat \beta)$ is never tangential to $\sum^{(1)}_{tan}$.
\end{observation}

Let $l$ be an integer and let us consider $l$ iterations of the loop $rot(\hat \beta)$. It follows from this lemma that the connected components of the natural resolution of $TL(l .rot(\hat \beta))$ (i.e. the abstract 
union of circles where the branches in the triple points are separated) are isotopy invariants of $\hat \beta \hookrightarrow V$.
We apply now our theory of one-cocycle polynomials but only to those triple points in $TL(l.rot(\hat \beta))$ where three 
{\em given \/} components (of the resolution) of $TL(l.rot(\hat \beta))$ intersect. The resulting invariants are called {\em character invariants}. Character invariants take into account the permutation (called {\em monodromy}) of the crossings of $\hat \beta$ which is induced by the loop $l.rot(\hat \beta)$. A character invariant for a fixed $l$ corresponds now to an unordered set of one-cocycle polynomials, namely one polynomial for each fixed triple of components (not necessarily different) of the resolution of $TL(l.rot(\hat \beta))$. The character invariants for $l=1$ coincide with the one-cocycle polynomials, as follows from Lemma 2 (see Section 2.3). Taking $l$ multiples of the loop $rot(\hat \beta)$ changes the character invariants in an uncontrollable way because of the monodromy of the crossings, which depends strongly on the closed braid $\hat \beta$. Hence, each character invariant  forms no longer a one-cocycle, but just a map $\mathbb{Z} \rightarrow$ {\em \{unordered sets of integer Laurent polynomials\}}.
 These new invariants are not of finite type but they are still calculable with polynomial complexity with respect to the braid length.

Alexander Stoimenow has written a computer program in order to calculate the simplest character invariants (see Section 4.3). It turns out that already these character invariants of linear complexity can sometimes detect the non-invertibility of closed braids  (i.e. the closed braid together with the core of the complementary solid torus in $S^3$ is a non-invertible link in $S^3$). 

The basic notions of our one parameter approach to knot theory, namely the space of non-singular knots, its 
discriminant, the stratification of the discriminant, the unfoldings of the strata in terms of singularity theory, the coorientation of strata of low codimension, the canonical loop, the trace graph, the equivalence relation for trace graphs,  are worked out in all details in our joint work with Vitaliy Kurlin \cite{F-K}, \cite{F-K2}. Therefore we concentrate in this paper only on the construction of the new invariants.

\vspace{0,5cm}

{\em Acknowledgments \/}

{\em I am grateful to Stepan Orevkov and Vitaliy Kurlin for many interesting discussions and to Allen Hatcher for a proof of Proposition 3. I am especially grateful to Alexander Stoimenow for writing his computer program and for calculating interesting examples.}

\section{Basic notions of one parameter knot theory}
In this section we recall briefly the basic notions of our theory. All details with complete proofs can be found in 
\cite{F-K} and \cite{F-K2} (even in the much more general setting of knots in 3-space).

\subsection{The space of diagrams of closed braids and its discriminant }

We work in the smooth category and all orientable manifolds are actually oriented.
We fix once for all a coordinate system in $\mathbb{R}^3$ : $(\phi , \rho , z )$.
Here, $(\phi ,\rho ) \in S^1\times \mathbb{R}^+$ are polar coordinates of the plane  $\mathbb{R}^2 = \{ z = 0 \}$.
A closed n-braid $\hat \beta$ is a knot in the solid torus $V = \mathbb{R}^3 \setminus  z-axes$, such that 
$\phi : \hat \beta \to  S^1$ is non-singular and $[\hat \beta] = n \in H_1(V)$. Let $M(\hat \beta)$ be the infinite dimensional 
space of all closed braids (with respect to $\phi$) which are isotopic to $\hat \beta$
in V. Let $M_n$ be the disjoint union of all spaces $M(\hat \beta)$.
A well known theorem of Artin (see e.g. \cite{M78}) says that two closed braids in the solid torus are isotopic as links 
in the solid torus 
if and only if they are isotopic as closed braids. Therefore it is enough to consider only isotopies through closed braids.
Let $pr : \mathbb{R}^3 \setminus z-axes \to \mathbb{R}^2 \setminus 0$ be the canonical projection 
$(\phi ,\rho , z) \to (\phi ,\rho )$.Each closed braid is then represented by a knot diagram with respect to $pr$. A generic 
closed braid $\hat \beta$ has only ordinary double points 
as singularities of $pr(\hat \beta)$. 

The discriminant $\sum$ has a natural stratification:  $\sum = \sum^{(1)} \cup \sum^{(2)} \cup ...$,
where $\sum^{(i)}$ are the union of all strata of codimension i in $M(\hat \beta)$.

\begin{theorem}
(Reidemeisters theorem for closed braids)

$\sum^{(1)} = \sum^{(1)}_{tri} \cup \sum^{(1)}_{tan}$,

where $\sum^{(1)}_{tri}$ is the union of all strata which correspond to diagrams with exactly one ordinary triple point 
(besides ordinary double points)  and $\sum^{(1)}_{tan}$ is the union of all strata which correspond to diagrams with exactly 
one 
ordinary autotangency.
\end{theorem}

In the sequel we need also the description of $\sum^{(2)}$.

\begin{theorem}
(Higher order Reidemeister theorem for closed braids)

$\sum^{(2)} = \sum ^{(2)}_{quad} \cup \sum ^{(2)}_{trans-self} \cup \sum ^{(2)}_{self-flex} \cup \sum^{(2)}_{inter}$ 

where $\Sigma^{(2)}_{quad}$ is the union of all strata which correspond to diagrams with exactly one  ordinary quadruple point,
$\Sigma^{(2)}_{trans-self}$ is the union of all strata which correspond to diagrams with exactly one ordinary autotangency 
through which passes another branch transversely,
$\Sigma^{(2)}_{self-flex}$ corresponds to the union of all strata of diagrams with an autotangency in an ordinary flex,
$\sum^{(2)}_{inter}$ is the union of all transverse intersections of strata from $\sum^{(1)}$.
\end{theorem}

\begin{figure}
\centering
\includegraphics{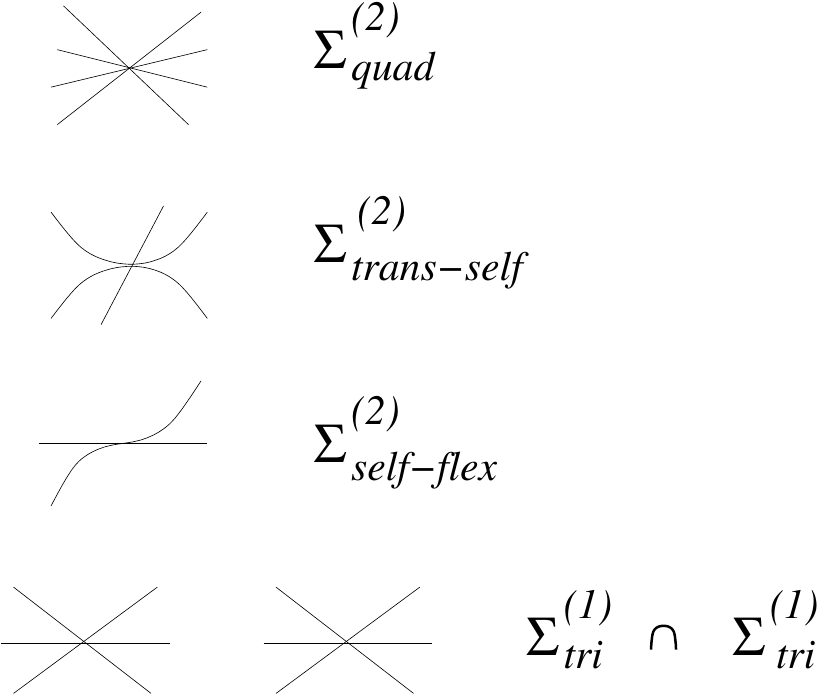}
\caption{\label{sing} The strata of codimension 2}  
\end{figure}

We show these strata in Fig.~\ref{sing}.

For a proof as well as for all other necessary preparations from singularity theory see \cite{F-K} (and also \cite{CS} and references therein).

Our strategy is the following: for an oriented  generic loop or arc in $M_n$ we associate some polynomial to the intersection with each stratum in $\Sigma^{(1)}_{tri}$, i.e. to each Reidemeister move of type III, and we sum up over all Reidemeister moves in the arc. It could be seen as an integration of a discrete closed 1-form over a discrete loop. We have to prove now that this sum is 0 for each meridian of strata in $\Sigma^{(2)}$. It follows that our sum is invariant under generic homotopies of arcs (with fixed endpoints). But it takes its values in an abelian ring and hence it is a 1-cocycle.

\subsection{The canonical loop}

We identify $\mathbb{R}^3 \setminus z-axes$ with the standard solid torus $V = S^1 \times D^2 \hookrightarrow \mathbb{R}^3 \setminus z-axes$.
We identify the core of V with the unit circle in $\mathbb{R}^2$.

Let $rot(V)$ denote the $S^1$-parameter family of diffeomorphismes of V which is defined in the following way:
we rotate the solid torus monotonously and with constant speed around its core by the angle t , $t \in [0 ,2\pi]$, i.e.
all discs $( \phi = const ) \times D^2$ stay invariant and are rotated simultaneously around their center.

 Let $\hat \beta$ be a closed braid.

\begin{definition}
The {\em canonical loop \/}  $rot(\hat \beta) \in M(\hat \beta)$ is the oriented loop induced by $rot(V)$.
\end{definition}

Notice that the whole loop $rot(\hat \beta)$ is completely determined by an arbitrary point in it.

The following lemma (see \cite{F-K}) is an immediate corollary of the definition of the canonical loop .
\begin{lemma} Let $\hat \beta_s , s\in [0 ,1 ]$ , be an isotopy of closed braids in the solid torus. Then $rot(\hat \beta_s ) , s \in [0 ,1 ]$,
is a homotopy of loops in $M(\hat \beta )$.
\end{lemma}
Evidently, the canonical loop can be defined for an arbitrary link in V in exactly the same way. However, in the case of closed braids
we can give an another (combinatorial) definition, which makes concrete calculations much easier.

Let $\Delta \in B_n$ be Garside's element, i.e. 

$\Delta=( \sigma_1\sigma_2\ ..\sigma_{n-1}) ( \sigma_1\sigma_2\ ..\sigma_{n-2}).. (\sigma_1 \sigma_2) (\sigma_1)$. Its square $\Delta^2$ is a generator of the center of $B_n$ (see e.g. \cite{B74}).
Geometrically, $\Delta^2$ is the full twist of the n strings.
\begin{definition}
Let $\gamma \in B_n$ be a braid with closure isotopic to $\hat \beta$. Then the {\em combinatorial canonical loop \/}
$rot(\gamma)$ is defined by the following sequence of braids:

$\gamma \to \Delta\Delta^{-1}\gamma \to \Delta^{-1}\gamma\Delta \to \dots \to \Delta^{-1}\Delta\gamma' \to \gamma' 
\to \Delta\Delta^{-1}\gamma' \to \Delta^{-1}\gamma'\Delta \to \dots \to \Delta^{-1}\Delta\gamma \to 
\gamma$

Here, the first arrow consists only of Reidemeister II moves, the second arrow is a cyclic permutation of the braid word 
(which corresponds to an isotopy of the braid diagram in the solid torus) and the following arrows consist of "pushing $\Delta$
monotonously from the right to the left through the braid $\gamma$". We obtain a braid $ \gamma'$, which is just $\gamma$ with each generator $\sigma_i$ replaced by $\sigma_{n-i}$, and we start again.
\end{definition}
We give below a precise definition in the case $n = 3$. The general case is a straightforward generalization which is left to the reader.
$\Delta = \sigma_1\sigma_2\sigma_1$ for $n = 3$. We have just to consider the following four cases:

$\sigma_1\Delta = \sigma_1(\sigma_1\sigma_2\sigma_1) \to \sigma_1(\sigma_2\sigma_1\sigma_2) = \Delta\sigma_2$

$\sigma_2\Delta = (\sigma_2\sigma_1\sigma_2)\sigma_1 \to (\sigma_1\sigma_2\sigma_1)\sigma_1 = \Delta\sigma_1$

$\sigma_1^{-1}\Delta = \sigma_1^{-1}(\sigma_1\sigma_2\sigma_1) \to \sigma_2\sigma_1 \to (\sigma_1\sigma_1^{-1})\sigma_2\sigma_1 \to \sigma_1(\sigma_2\sigma_1\sigma_2^{-1}) = \Delta\sigma_2^{-1}$

$\sigma_2^{-1}\Delta = \sigma_2^{-1}(\sigma_1\sigma_2\sigma_1) \to (\sigma_1\sigma_2\sigma_1^{-1})\sigma_1 \to \sigma_1\sigma_2 \to \sigma_1\sigma_2\sigma_1\sigma_1^{-1} = \Delta\sigma_1^{-1}$.

Notice, that the sequence is canonical in the case of a generator and almost canonical in the case of an inverse generator. Indeed,
we could replace the above sequence   $\sigma_1^{-1}\Delta \to \Delta\sigma_2^{-1}$ by

$\sigma_1^{-1}(\sigma_1\sigma_2\sigma_1) \to \sigma_2\sigma_1 \to \sigma_2\sigma_1\sigma_2\sigma_2^{-1} \to (\sigma_1\sigma_2\sigma_1)\sigma_2^{-1}$.

But it turns out that the corresponding canonical loops in $M(\hat \beta)$ differ just by a homotopy which passes once transversely through a stratum of $\sum^{(2)}_{trans-self}$ and our one-cocycles are invariant under this homotopy.

Let c be the word length of $\gamma$. Then we use exactly $2c(n-2)$ braid relations (or equivalently, Reidemeister III moves) in the combinatorial canonical loop.
This means that the corresponding loop in $M(\hat \beta)$ intersects $\sum^{(1)}_{tri}$ transversely in exactly $2c(n-2)$ points.

One easily sees that the combinatorial canonical loop $rot(\gamma)$ from Definition 2 differs from the geometrical 
canonical loop $rot(\hat \beta)$ from Definition 1 only by loops which correspond to rotations of the solid torus along its core (i.e. around the axis of the complementary solid torus in $S^3$).
But each such loop is just an isotopy of diagrams with respect to $pr$ and does not intersect the discriminant $\sum^{(i)}, i>0,$ at all.

\subsection{The trace graph}

The trace graph $TL(rot(\hat \beta))$ is our main combinatorial object, compare \cite {F-K} and \cite{F-K2}. It is an oriented singular link in a thickened torus. All its singularities are ordinary triple points.

Let $\hat \beta_t , t \in S^1$, be the (oriented) family of closed braids corresponding to the canonical loop $rot(\hat \beta)$.
We assume that the loop $rot(\hat \beta)$ is a generic loop. Let $\{ p_1^{(t)}, p_2^{(t}, \dots , p_m^{(t)} \}$ be the set of double points of
$pr(\hat \beta_t) \subset S_\phi^1 \times \mathbb{R}_\rho$.
The union of all these crossings for all $t \in S^1$ forms a link $TL(rot(\hat \beta)) \subset (S_\phi^1 \times \mathbb{R}_\rho^+) \times S_t^1$ (i.e. we forget the coordinate $z(p_i^{(t)})$).
$TL(rot(\hat \beta))$ is non-singular besides ordinary triple points which correspond exactly to the triple points in the family $pr(\hat \beta_t)$.
A generic point of $TL(rot(\hat \beta))$ corresponds just to an ordinary crossing $p_i^{(t)}$ of some closed braid $\hat \beta_t$.
Let $t: TL(rot(\hat \beta)) \to S_t^1$ be the natural projection. We orient the set of all generic points in $TL(rot(\hat \beta))$
(which is a disjoint union of embedded arcs) in such a way that the local mapping degree of t at $p_i^{(t)}$ is $+1$ if and only if $p_i^{(t)}$
is a positive crossing (i.e. it corresponds to a generator of $B_n$ , or equivalently , its {\em writhe \/} $w(p_i^{(t)}) = +1$).

The arcs of generic points come together in the triple points and in points corresponding to an ordinary autotangency in some $pr(\hat \beta)$.
But one easily sees that the above defined orientations fit together to define an orientation on the natural resolution $TL\tilde(rot (\hat \beta))$ of $TL(rot(\hat \beta))$, compare also \cite{F01}. $TL\tilde(rot(\hat \beta))$ is a union of oriented circles, called {\em trace circles \/}.We can attach {\em stickers \/} $i \in \{ 1, 2, \dots ,n-1 \}$ to the edges of $TL(rot(\hat \beta))$ in the following way:
each edge of $TL(rot(\hat \beta))$ corresponds to a letter in a braid word. Indeed, each generic point in an edge corresponds to an ordinary crossing of a braid 
projection and , hence , to some $\sigma_i$ or some $\sigma_i^{-1}$ . We attache to this edge the number i . The information about the exponent $+1$
or $-1$ is contained in the orientation of the edge.

We identify $H_1(V)$ with $\mathbb{Z}$ by sending the homology class which is generated by the closed 1-braid to the generator $+1$. If $\hat \beta$ is a knot then we can attache to each trace circle a 
{\em homological marking \/} $a \in H_1(V)$ in the following way: let p be a crossing corresponding to a generic point in the trace circle.
We smooth p with respect to the orientation of the closed braid. The result is an oriented 2-component link. The component of this link which contains 
the under-cross  which goes to the over-cross at p is called $p^+$. We associate now to p the homology class $a = [p^+] \in H_1(V)$, compare also \cite{F93}.
One easily sees that $a \in \{ 1 ,2 , \dots , n-1\}$ and that the class a does not depend on the choice of the generic point in the trace circle.
Indeed, the two crossings involved in a Reidemeister II move have the same homological marking and a Reidemeister III move 
does not change the homological marking of any of the three involved crossings. 

We could see $TL(rot(\hat \beta))$ as an object which takes into account simultaneously all possible directions of projections of $\hat \beta$ into isotopic properly embedded annuli in $V$ which contain the core of $V$.

Evidently, crossings with different homological markings belong to different trace circles. Surprisingly, the inverse is also true, see  Lemma 3.1 in \cite{F-K2}.

\begin{lemma}
Let $TL(rot(\hat \beta))$ be the trace graph of the closure of a braid $\beta \in B_n$, such that $\hat \beta$ is a knot.
Then $TL\tilde (rot(\hat \beta))$ splits into exactly n-1 trace circles. They have pairwise different homological markings.
\end{lemma}
Consequently, the trace circles are characterized by their homological markings. Notice, that the set of homological markings is independent of the word length of the braid (which has a knot as closure).

\subsection{A higher order Reidemeister theorem for trace graphs of closed braids }

\begin{figure}
\centering
\includegraphics{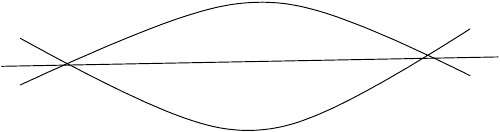}
\caption{\label{trihed} A trihedron}  
\end{figure}

\begin{definition}
A {\em trihedron \/} is a 1-dimensional subcomplex of $TL(rot(\hat \beta))$ which is contractible in the thickened torus and which has the form as shown in Fig.~\ref{trihed}.
\end{definition}

\begin{figure}
\centering
\includegraphics{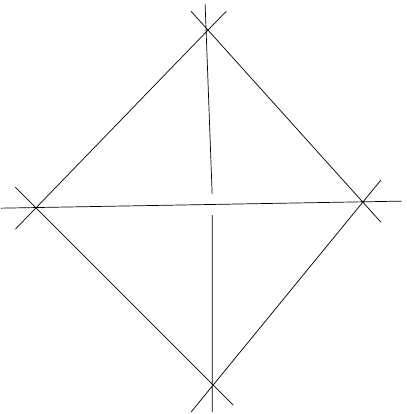}
\caption{\label{tetra} A tetrahedron}  
\end{figure}

\begin{definition}
A {\em tetrahedron \/} is a 1-dimensional subcomplex of $TL(rot(\hat \beta))$ which is contractible in the thickened torus and which has the form as shown in Fig.~\ref{tetra}.

\end{definition}

\begin{figure}
\centering
\includegraphics{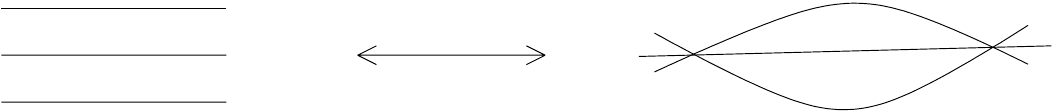}
\caption{\label{trimove} A trihedron move}  
\end{figure}

\begin{definition}A {\em trihedron move \/} is shown in Fig.~\ref{trimove}.
\end{definition}

\begin{figure}
\centering
\includegraphics{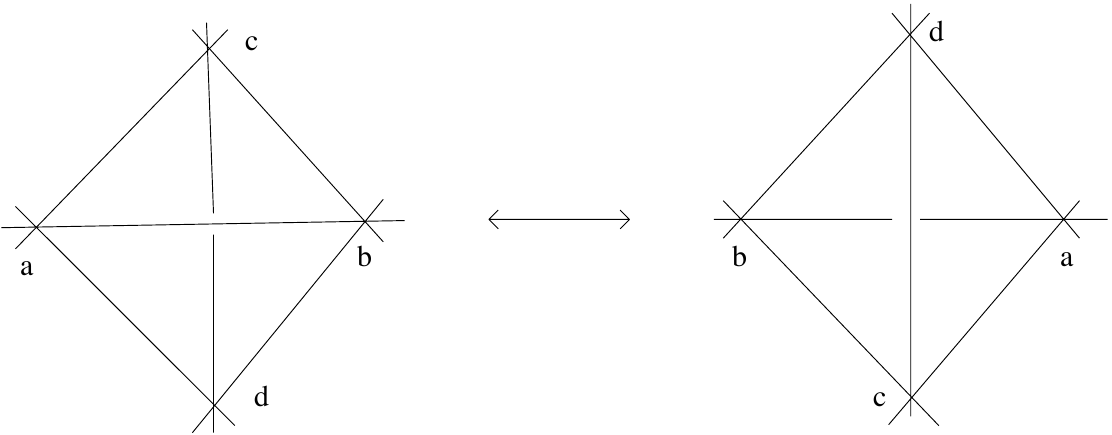}
\caption{\label{tetramove} A tetrahedron move}  
\end{figure}

\begin{definition}
A {\em tetrahedron move \/} is shown in Fig.~\ref{tetramove}.
\end{definition}

The rest of $TL(rot(\hat \beta)) \hookrightarrow S^1 \times S^1 \times \mathbb{R}^+$ is unchanged under the moves. The stickers on the edges change in the canonical way.

Notice, that a trihedron move corresponds to a generic homotopy of the canonical loop which passes once through an ordinary tangency with a stratum of $\sum^{(1)}_{tri}$. A tetrahedron move corresponds to a generic homotopy of the canonical loop which passes transversely once through a stratum of $\sum^{(2)}_{quad}$, i.e. corresponding to an ordinary quadruple point in the projection.

\begin{definition}
The {\em equivalence relation for trace graphs $TL(rot(\hat \beta))$ \/} is generated by the following three operations:

(1) isotopy in the thickened torus

(2) trihedron moves

(3) tetrahedron moves 
\end{definition}

The following important Reidemeister type theorem is a particular case of Theorem 1.10.  in \cite{F-K}.

\begin{theorem}
(Higher order Reidemeister theorem for trace graphs of closed braids)
 
Two closed braids (which are knots) are isotopic in the solid torus if and only if their trace graphs in the thickened torus are equivalent.
\end{theorem}

\begin{remark}
Notice, that not all representatives of an equivalence class of a trace graph correspond to the canonical loop (which is very rigid, because it is determined by a single braid diagram in it) of some closed braid. 
However, one easily sees that all representatives correspond to loops in the space $M(\hat \beta)$.
\end{remark}
\begin{remark}
A generic homotopy of loops $\gamma_s , s \in [0 , 1]$, in $M(\hat \beta)$ could of course become tangential for some $s$ to $\sum^{(1)}_{tan}$ at a generic point. We can allow tangencies of $\sum^{(1)}_{tan}$ on the negative side, i.e. where the ordinary diagrams have two crossings less, because this corresponds to the birth or the death of a (null-homologous) component of the trace graph. However, a tangency on the positive side would imply a Morse modification of index 1 of the trace graph and, hence, change the components of its natural resolution in an uncontrollable way. The important point is, that this does not happen for the (very rigid) homotopies of $rot(\hat \beta)$ which are induced by generic isotopies of $\hat \beta$ in V, compare Observation 1. 
\end{remark}

{\em Proof of Observation 1}. Indeed, under a monotonous rotation of the closed braid around the core of the solid torus $V$, $rot(\hat \beta)$ is tangential to $\sum^{(1)}_{tan}$ in an ordinary point if and only if the trace graph $TL(rot(\hat \beta))$  has a Morse singularity. This could only happen if two strings of the braid would be become tangential to the same disc 
$( \phi = const ) \times D^2$ in the fibration of $V$. But this is not possible, because the closed braids are always transverse to the disc-fibration of $V$.
 $\Box$

Consequently, the trace circles of (the resolution of) the trace graph are isotopy invariants for closed braids.

\begin{remark}
There can be loops $\gamma_s$ in a homotopy which are tangential to $\sum^{(1)}_{tri}$. For the trace graphs associated to the loops $\gamma_s$ this corresponds to a trihedron move (compare \cite {F-K} and \cite{F-K2}).
\end{remark}

We have some more information about isotopies of trace graphs.

\begin{figure}
\centering
\includegraphics{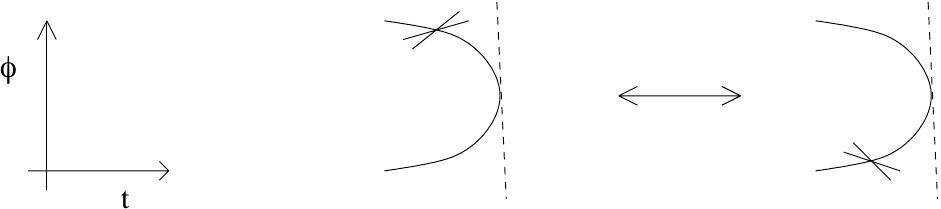}
\caption{\label{triptan} A triple point slides over a vertical tangency}  
\end{figure}

\begin{definition}
A {\em time section \/} in the thickened torus $(S^1_\phi \times \mathbb{R}^+_\rho) \times S^1_t$ is an annulus of the form 
$(S^1_\phi \times \mathbb{R}^+_\rho) \times \{ t = const \}$.
\end{definition}
The intersection of $TL(rot(\hat \beta))$ with a generic time section corresponds to the crossings of the closed braid $\hat \beta_t$.
Using the orientation and the stickers on $TL(rot(\hat \beta))$ we can read off a cyclic braid word for $\hat \beta$ in each generic time section.
The tangent points of $TL(rot(\hat \beta))$  with  time sections correspond exactly to the Reidemeister II moves in the one parameter family of diagrams $\hat \beta_t , t \in S^1$.
A triple point in $TL(rot(\hat \beta))$ slides over such a tangent point if and only if the canonical loop passes in a homotopy transversely through a stratum of $\sum^{(2)}_{trans-self}$.
We illustrate this in Fig.~\ref{triptan}.

\begin{figure}
\centering
\includegraphics{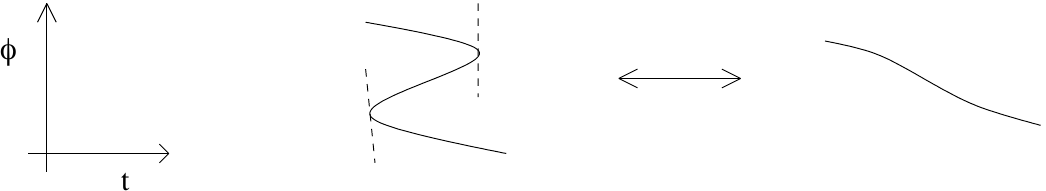}
\caption{\label{zigzag} A couple of tangencies appears or disappears}  
\end{figure}

When the canonical loop passes transversely through a stratum of $\sum^{(2)}_{self-flex}$ then the trace graph changes as shown in Fig.~\ref{zigzag}.

\begin{figure}
\centering
\includegraphics{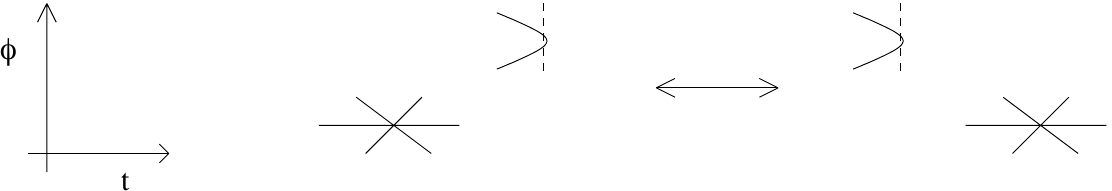}
\caption{\label{comm} A commutation of a triple point with a vertical tangency}  
\end{figure}

Finally, when the canonical loop passes transversely through a stratum of $\sum^{(2)}_{inter}$ then the t-values of triple points or tangencies with time sections are interchanged. We show an example in Fig.~\ref{comm}.

\vspace{0.4cm}
It turns out that we could replace Observation 1 by the following observation of Orevkov, but at the cost that braids become much longer.

Remember that in our one-parameter theory we consider (geometric) braids as tangles and not as isotopy classes of tangles (or elements in the braid group). So, a generic braid is for us a word in the standard generators and their inverses of the braid group.

If a braid $\beta$ is in {\em Garsides normal form}, see \cite{G}, then it is in particular of the form 
$\beta^+ (\Delta^2)^k$, where $\beta^+$ is a positive braid, $k$ is an integer and $\beta^+$ is not divisible by $\Delta^2$. 

Let $\beta^+ (\Delta^2)^k$ and $\beta^{'+} (\Delta^2)^k$ be two n-braids in Garsides normal form. {\em Garsides solution of the conjugacy problem} says that the the isotopy classes of the braids are conjugate in $B_n$  if and only if $\beta^+$ and $\beta^{'+}$ can be related by a sequence of the following operations (see \cite{G}):

1) isotopy amongst positive braids

2) cyclic permutation

3) conjugacy with permutation braids $s\beta^+s^{-1}$ followed by a simplification to a positive braid
(where the {\em permutation braids $s$} are exactly the left and right divisors in the braid monoid $B_n^+$
of Garsides braid $\Delta$, i.e. there exist positive n-braids $s'$ and $s''$ such that $ss'$ is isotopic to $s''s$ which is isotopic to $\Delta$).

Stepan Orevkov has made the following observation.

\begin{observation}(Orevkov)

Assume that the closure $\hat \beta^+$ contains a half-twist $\Delta$, i.e. $\beta^+$ is equivalent to $\beta^{''+}\Delta$  for some positive braid $\beta^{''+}$ by operations 1 and 2 only. Then we can skip the operation 3 in Garsides solution of the conjugacy problem.

\end{observation}

{\em Proof.} Indeed, assume that we have to perform the operation 3:

$\beta^+ \to s\beta^+s^{-1}=\beta^{'+}$.

From our hypothesis $\beta^+=\beta^{''+}\Delta$ by using operations 1 and 2. $\beta^{''+}\Delta=\beta^{''+}s^{''}s \to s\beta^{''+}s^{''}ss^{-1}=s\beta^{''+}s^{''} = \beta^{'+}$. But by using only operations 2 we can bring the last braid to $\beta^{''+}s^{''}s=\beta^{''+}\Delta=\beta^{+}$. Consequently, we have replaced the operation 3 by a sequence of only operations of type 1 and 2. $\Box$

Orevkov's observation implies that if two positive closed braids which contain a half-twist are isotopic in the solid torus then there is such an isotopy which stays within positive closed braids (because negative crossings occur only in operation 3). The loop $rot(\hat \beta^{''+}\Delta)$ can be represented by pushing twice $\Delta$ through $\beta^{''+}$. It follows easily now (by considering braids as diffeomorphisms of the punctured disc and $rot$ as the rotation  from 0 to $2\pi$ of the disc around its center) that if $\hat \beta^{''+}\Delta$ is isotopic to $\hat \beta^{'+}\Delta$ then the loops $rot(\hat \beta^{''+}\Delta)$ and $rot(\hat \beta^{'+}\Delta)$ are homotopic through positive closed braids (i.e. all generic crossings are positive). Consequently, there are no Morse modifications of the trace graph at all in such a homotopy and hence the trace circles are again isotopy invariants of the closed braid.

Adding full-twists to a braid does not change its geometry and its geometric invariants (as hyperbolic volume or entropy for pseudo-Anosov braids) and a full twist (which generates the center of the braid group) in a closed braid can be easily detected. Using Orevkov's observation we could therefore restrain our-self to the heart of the matter: constructing 1-cocycles for positive closed braids in $M(\hat \beta^+ \Delta)$, but of course braids would be in general much longer. Therefore we allow in this paper negative crossings as well.

\section{One-cocycle polynomials}

In this section we introduce our one-cocycle polynomials. Notice, that we do not need trace graphs here. They will only be needed in the refinements to character invariants.

\subsection{Gauss diagrams for closed braids with a triple crossing }
Let $f : S^1 \to \hat \beta$ be a generic orientation preserving diffeomorphisme. Let p be any crossing of $\hat \beta$.
We connect $f^{(-1)}(p) \in S^1$ by an oriented arrow, which goes from the under-cross to the over-cross and we decorate it by the writhe $w(p)$.
Moreover, we attache to the chord the homological marking $[p^+]$ (compare Section 2.3). The result is called a {\em Gauss diagram \/} for
$\hat \beta$ (compare also e.g. \cite{PV} and \cite{F01}).

One easily sees that $\hat \beta$ up to isotopy is determined by its Gauss diagram and the number $n = [\hat \beta] \in H_1(V)$.

We can form oriented loops in Gauss diagrams of closed braids, which are knots, in the following way: going along the circle following its orientation we can possibly jump at arrows and continue going along the circle following its orientation up to reaching our starting point. The following simple observation is at the origin of our new solution of the tetrahedron equation which leads immediately to our polynomial valued one-cocycles.

\begin {observation}
Each such loop in the Gauss diagram of a closed braid, which is a knot, represents a homology class in $\{ 1 ,2 , \dots , n-1\}$.
\end {observation}

{\em Proof.} Each such loop is positive transverse to the disc fibration of $V$ exactly as the closed braid, besides at the jumps, which can be represented by arcs in the discs. Consequently, the loop represents a homology class between $1$ and $n-1$.
$\Box$

A {\em Gauss sum of degree d \/} is an expression assigned to a diagram of a closed braid which is of the following form, compare \cite{F01}:

$\sum$ function( writhes of the crossings)

where the sum is taken over all possible choices of d (unordered) different crossings in the knot diagram such that the arrows without the writhes arising from these crossings 
build a given sub-diagram with given homological markings. The marked sub-diagrams (without the writhes) are called {\em configurations\/}.
If the function is the product of the writhes (this is always the case in this paper), then we will denote the sum shortly by the configuration itself.
We need to define Gauss diagrams for knots with an ordinary triple point too. The triple point corresponds to a triangle in the Gauss diagram of the knot.
Notice, that the preimage of a triple point has a natural ordering coming from the orientation of the $\mathbb{R}^+$-factor.
One easily sees that this order is completely determined by the arrows in the triangle.

We provide each stratum of $\sum^{(1)}_{tri}$ with a co-orientation which depends only on the non-oriented  underlying curves $pr(\hat \beta)$ in $S^1 \times \mathbb{R}^+$.
Consequently, for the definition of the co-orientation we can replace the arrows in the Gauss diagram simply by chords.

\begin{figure}
\centering
\includegraphics{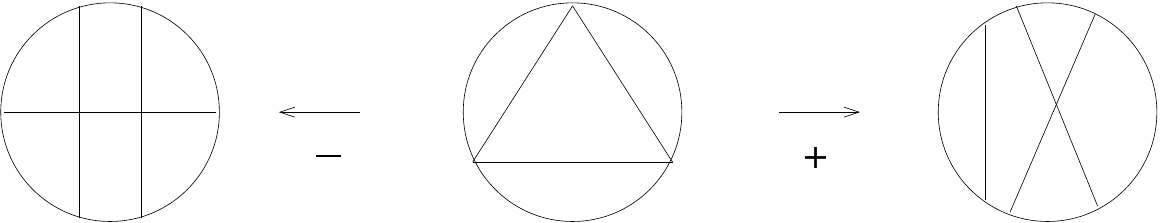}
\caption{\label{coorient} The coorientation for a triple crossing}  
\end{figure}

\begin{definition}
The {\em co-orientation \/} of the strata in $\sum^{(1)}_{tri}$ is given in Fig.~\ref{coorient}. 
\end{definition}

Let $s$ be a point in the transverse intersection of an oriented arc $S \subset M_n$ with $\sum^{(1)}_{tri}$. Then the sign $sign(s)$ is $+1$ if the orientation of $S$ agrees with the co-orientation of $\sum^{(1)}_{tri}$ at $s$ and  $sign(s)=-1$ otherwise.

\begin{figure}
\centering
\includegraphics{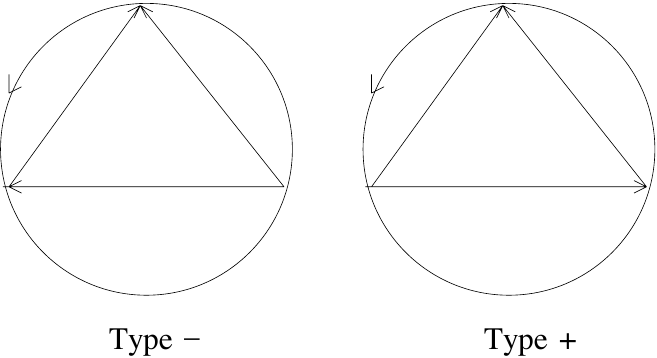}
\caption{\label{globtype} The unmarked global types of triple crossings}  
\end{figure}

There are exactly two types of triple points without markings. We show them in Fig.~\ref{globtype}.

\begin{figure}
\centering
\includegraphics{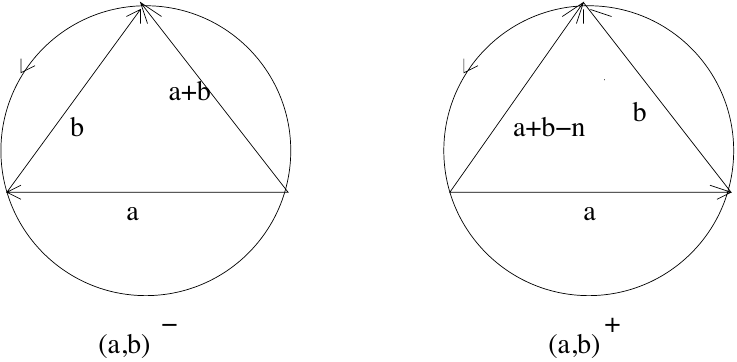}
\caption{\label{marktype} The marked global types of triple crossings}  
\end{figure}

We attach now the homological markings to the three chords. Let $a,b \in \{1 , 2 , \dots , n-1\}$ be fixed. 
Then the markings of a triple crossing are as shown in Fig.~\ref{marktype}. We encode the types of the marked triple points by $(a ,b)^-$ and respectively $(a ,b)^+$.
The union of the corresponding strata of $\sum^{(1)}_{tri}$ are encoded in the same way.

\subsection{The general form of one-cocycle polynomials from Gauss diagrams with a triple crossing}

We will construct one-cocycles on the space $M_n$ (the space of all closed n-braids which are knots).
We obtain invariants of closed braids when we evaluate these cocycles on the homology class represented by the canonical loop $rot$. (If a braid $\beta$ is reducible then there is an incompressible not-boundary parallel  torus in the complement $V \setminus \hat \beta$. This implies that there are new loops in $M(\hat \beta)$ to which we could apply our one-cocycles as well, compare \cite{H} for the case of long knots.)
 
Let $n, d \in \mathbb{N}^*$ be fixed. Let $(a,b)^+$, respectively $(a,b)^-$, be a fixed
{\em type of marked triple point\/} as shown in Fig.~\ref{marktype}.

\begin{definition}
A {\em configuration I of degree d\/} is an abstract Gauss diagram
without writhes which contains exactly $d$ arrows marked in $\{ 1 ,2 , \dots , n-1\}$ besides the triangle $(a,b)^+$ or $(a,b)^-$.
\end{definition}

Let $\{ I_i \}$ be the finite set of all configurations of degree $d$
with respect to $(a,b)^\pm$. Let $\Gamma = \sum_{i}
{\epsilon_i I_i}$ be a given linear combination with each $\epsilon_i \in
\{ 0, +1, -1 \}$. (The type of the triple point is {\em always\/}
fixed in any cochain $\Gamma$).

Here is the general form of our cochains.

\begin{definition}
$\Gamma$ gives rise to a {\em 1-cochain\/}
by assigning to each oriented generic loop $S \subset M_n$
an integer Laurent polynomial $\Gamma(S)$ in the following way:
\begin{displaymath}
\Gamma(S) =
\sum_{s_i \in S \cap \Sigma^{(1)}_{tri} \textrm{of type $(a,b)^\pm$}}
 {sign(s_i)x^{\sum_{i}{\epsilon_i (\sum_{D_i}{\prod_{j}{w(p_j)})}}}}
\end{displaymath}
where $D_i$ is the set of unordered $d$-tuples $(p_1, \dots, p_d)$
of arrows which enter in $I_i$ in the Gauss diagram of $s_i$ .
\end{definition}

\begin{lemma}
 If  $\Gamma(S)$ is invariant under each generic deformation
of $S$ through any stratum of $\Sigma^{(2)}$, then $\Gamma$
is a 1-cocycle.
\end{lemma}
{\em Proof:\/}  In this case, $\Gamma(S)$ is
invariant under homotopies of $S$. Indeed, tangent points of $S$ with $\Sigma^{(1)}_{tri}$
correspond just to trihedron moves. The two triple points give the same contribution to 
$\Gamma(S)$ but with different signs. A tangency with $\Sigma^{(1)}_{tan}$ does not change the contribution
of the triple points at all. This implies the invariance under
homologies of $S$  because $\Gamma(S)$ takes it values in an abelian ring. $\Box$

\begin{definition}
A {\em cohomology class\/} in $H^1(M_n; \mathbb{Z}[x,x^{-1}])$ is {\em
of Gauss degree $d$\/} if it can be represented by some 1-cocycle
$\Gamma$ such that each configuration $I_i$ has at most  $d$ arrows (besides the three arrows of the triangle).
\end{definition}

\begin{remark}

The integer valued invariant $d(\Gamma(rot(\hat\beta))/dx$ at $x=1$ is given by a Gauss diagram formula and it is hence an invariant of finite type. Indeed, the one-cocycle invariant
$\Gamma(rot(\hat\beta))$ is calculated as some sum $\sum_{s_i}$ over
triple points $s_i$ in $rot(\hat\beta)$. Therefore, it suffices to
prove that this sum $\sum_{s_i}$ for each triple point $s_i$ is of
finite type (even if it is not invariant).
If $\Gamma$ is of Gauss degree $d$ then $\sum_{s_i}$ depends only on the triple
point and of configurations of $d$ other crossings. This means that
in order to calculate a summand in $\sum_{s_i}$ we can switch all other
crossings besides the triple point and the fixed $d$ crossings. The
result will not change. This implies immediately that each $\sum_{s_i}$
is of degree at most $d+1$ (see \cite{PV} and also \cite{F01}).\end{remark}

The above definition induces a filtration on a part of $H^0(M_n;\mathbb{Z})$
by taking all $d(\Gamma(rot(\hat\beta)))/dx \vert_{x=1}$.

Let $M$ be the (disconnected) space of all embeddings $f: S^1
\hookrightarrow \mathbb{R}^3$. Vassiliev \cite{V90} has introduced a
filtration on a part of $H^0(M; \mathbb{Z})$ using the discriminant
$\Sigma_{sing}$ of singular maps. It is not difficult to see
that each component of the space of all (unparametrized) differentiable maps of the
circle into the solid torus is contractible. Indeed, there is an
obvious canonical homotopy of each (perhaps singular) knot to a
multiple of the core of the solid torus. The core of the solid torus is invariant under $rot_{S^1}(V)
\times rot_{D^2}(V)$. Thus, the above space is star-like. Therefore,
Alexander duality could be applied and Vassilievs original approach could be
generalized for knots in the solid torus too.
It would be interesting to compare his filtration with our
filtration (we compare them just for the Gauss degree 0 in Section 3.3).

In the next sections, we will construct 1-cocycle polynomials $\Gamma$
in an explicit way. We show with a simple example that they are in general not of finite type.

\subsection{One-cocycles of Gauss degree 0}

Let $\beta \in B_n$ be such that its closure $\hat \beta \hookrightarrow
V$ is a knot.

Let us recall the simplest finite type invariants for closed braids in Vassiliev's sense.

\begin{proposition}
 The space of finite type invariants of degree 1 is of dimension
$[n/2]$ (here $[.]$ is the integer part). It is generated by the
Gauss diagram invariants $W_a(\hat \beta) = \sum_{}{w(p)}$, where $a \in
\{ 1, 2, \dots, [n/2] \}$. The sum is over all crossings with fixed homological marking a.
\end{proposition}
{\em Proof:\/} It follows from Goryunov's \cite{G96} generalization of
finite type invariants for knots in the solid torus (likewise by generalizing the Kontsevich integral or by generalizing the defining skein relations for finite type invariants) that the invariants
of degree 1 correspond just to marked chord diagrams with only one
chord. Obviously, all these invariants can be expressed as Gauss
diagram invariants:

$W_a(\hat \beta) = \sum_{}{w(p)}$, $a \in \{ 1, \dots, n-1 \}$

\vspace{0.5cm}

(see also \cite{F93}, and Section 2.2 in \cite{F01}.)

Let us define $V_a(\hat \beta) := W_a(\hat \beta) - W_{n-a}(\hat \beta)$
for all $a \in \{ 1, \dots, n-1 \}$. We observe that $V_a(\hat \beta)$
is invariant under switching crossings of $\hat \beta$.
Indeed, if the marking of the crossing $p$ was $[p] = a$, then the
switched crossing $p^{-1}$ has marking $[p^{-1}] = n-a$, but $w(p)=
-w(p^{-1})$.
But every braid $\beta \in B_n$ is homotopic to $\gamma = \prod_{i=1}
^{n-1}{\sigma_i}$. A direct calculation for $\gamma$ shows that
$V_a(\hat \gamma) \equiv 0$. It is easily shown by examples that
$W_a, a \in \{1, \dots, [n/2] \}$ (seen as invariants in with values in $\mathbb{Q}$) are linearly independent. $\Box$

\begin{lemma}
 Let $a, b \in \{ 1 ,2 , \dots , n-1\}$ be fixed.
Consider the union of all co-oriented strata of $\Sigma^{(1)}$
which correspond to triple points of type either $(a, b)^-$ or $(a, b)^+$.
The closure in $M_n$ of each of these sets form 
integer cycles of codimension 1 in $M_n$.
\end{lemma}
\begin{remark}
Otherwise stated, $\Gamma_{(a,b)^+}$  and
$\Gamma_{(a,b)^-}$ both define integer valued 1-cocycles
of Gauss degree 0.
$\Gamma_{(a,b)^\pm}(S)$ is in this case by definition just the
algebraic intersection number of $S$ with the corresponding union
of strata of $\Sigma^{(1)}_{tri}$. The variable $x$ enters the one-cocycles only starting from the Gauss degree 1.
 
\end{remark}

\begin{figure}
\centering
\includegraphics{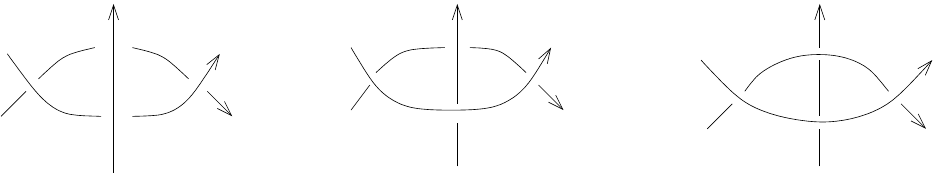}
\caption{\label{transself} Different triple crossings come together in an auto-tangency}  
\end{figure}

{\em Proof:\/} According to Section 2.1, we have to prove that the
co-oriented strata fit together in $\Sigma^{(2)}_{quad}$, $\Sigma^{(2)}_{trans-self}$ and $\Sigma^{(2)}_{inter}$. The first is evident, because at a stratum of $\Sigma^{(2)}_{quad}$ just four strata of $\Sigma^{(1)}_{tri}$ 
intersect pairwise transversely.
For the second, we have to distinguish 24 cases. Three of them are
illustrated in Fig.~\ref{transself}.

\begin{figure}
\centering
\includegraphics{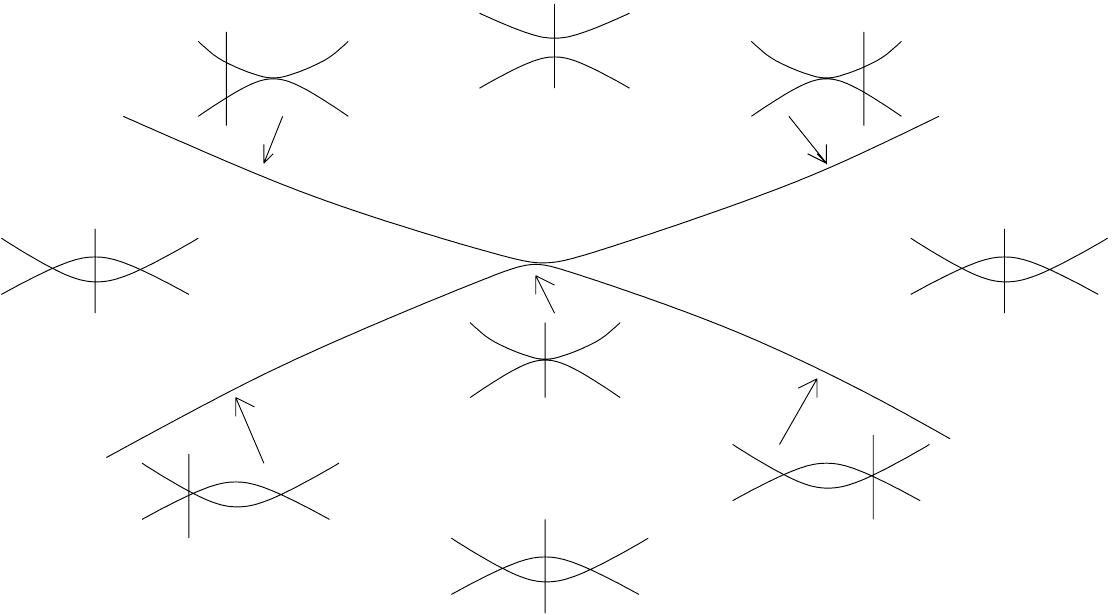}
\caption{\label{unftrans} The unfolding of an auto-tangency with a transverse branch}  
\end{figure}

The whole picture in a normal 2-disc of $\Sigma^{(2)}_{trans-self} \subset M_n$ is then shown in Fig.~\ref{unftrans}, but we draw only the corresponding planar curves.

All other cases are obtained from these three by inverting the orientation
of the vertical branch, by taking the mirror image (i.e. switching
all crossings), and by choosing one of two possible closings of the
3-tangle (in order to obtain an oriented knot).
In all cases, one easily sees that the two adjacent triple points
are always of the same marked type and that the co-orientations fit
together. $\Box$

\begin{proposition}
$\Gamma_{(a,b)^+}$  defines a non-trivial 1-cohomology
class of Gauss degree 0 if and only if $a \not= b$ and $a + b \leq n-1$.

$\Gamma_{(a,b)^-}$  defines a non-trivial 1-cohomology
class of Gauss degree 0 if and only if $a \not= b$ and $a + b \geq n+1$.

The following identities hold:

(*) \hspace{4cm}    $\Gamma_{(a,b)^+} + \Gamma_{(b,a)^+} \equiv 0$

 (*) \hspace{4cm}   $\Gamma_{(a,b)^-} + \Gamma_{(b,a)^-} \equiv 0$

\end{proposition}

{\em Proof:\/} For closed braids, the markings are all in
$\{ 1, \dots, n-1 \}$. Therefore, if
$a+b>n-1$ in $\Gamma_{(a,b)^+}$ or $a+b<n+1$ in $\Gamma_{(a,b)^-}$,
then there is no such triple point at all and the 1-cocycle is trivial. It follows from the identities that $\Gamma_{(a,a)^+}$
and $\Gamma_{(a,a)^-}$ are trivial.

Examples show that all the remaining 1-cocycles are non-trivial.

\begin{figure}
\centering
\includegraphics{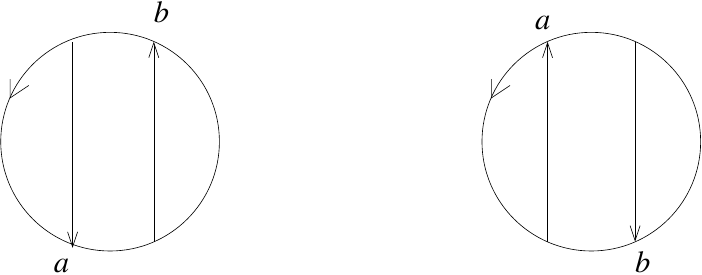}
\caption{\label{subdia} Two sub-diagrams}  
\end{figure}

In order to prove the identities, we use the following Gauss
diagram sums (see also section 1.6 in \cite{F01}):
\begin{displaymath}
I^+_{(a,b)} = \sum_{} w(p)w(q), \qquad I^-_{(a,b)} = \sum_{} w(p)w(q)
\end{displaymath}
Here, the first sum is over all couples of crossings which form a sub-configuration as shown on the left in Fig.~\ref{subdia}. The second sum is over all couples of crossings which form a sub-configuration as shown on the right in Fig.~\ref{subdia}. (Here, a and b are the homological markings.)

\begin{figure}
\centering
\includegraphics{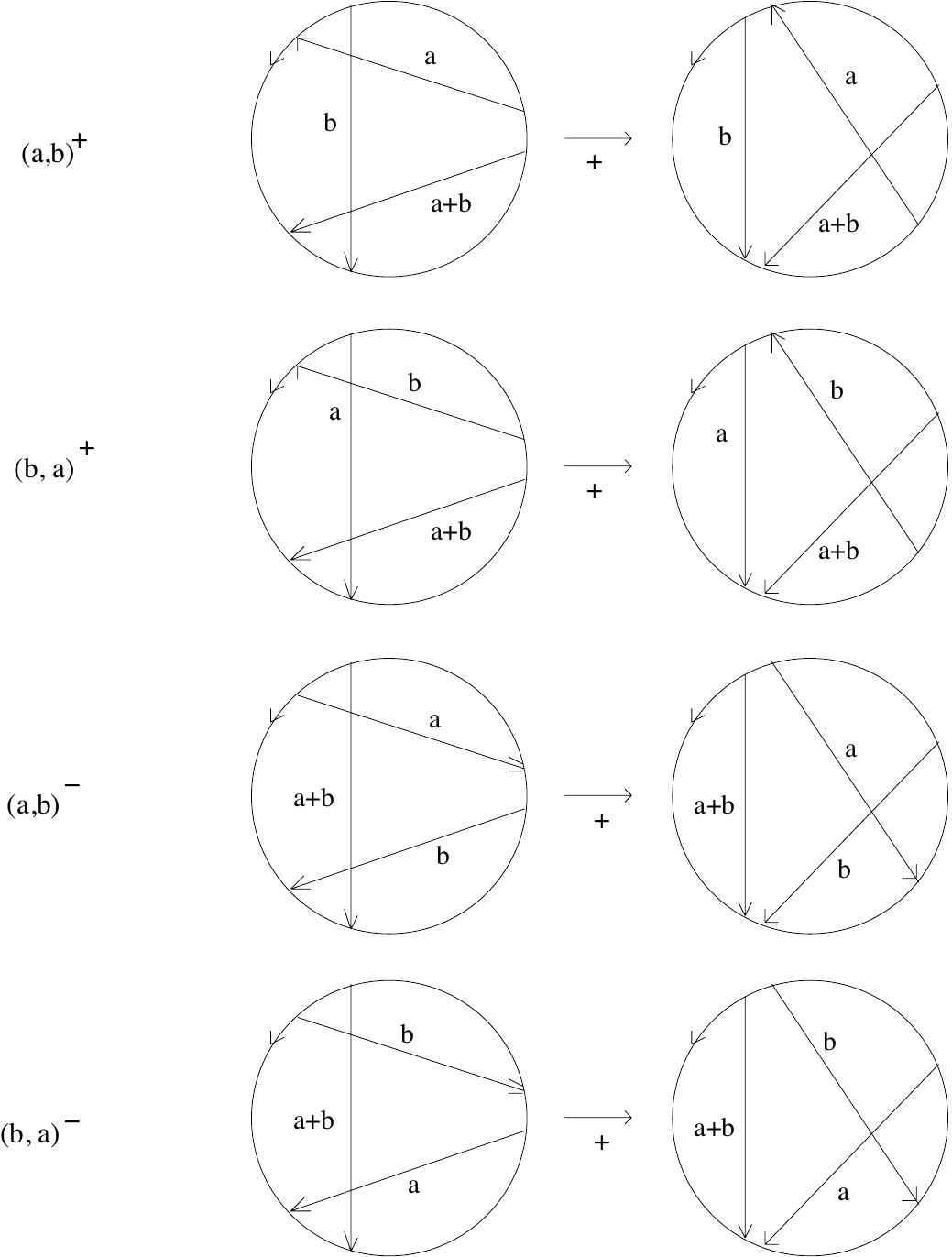}
\caption{\label{perturb} Perturbations of triple crossings}  
\end{figure}

These sums applied to diagrams of $\hat \beta$ are not invariants.
Let $S \subset M(\hat \beta)$ be a generic loop. Then $I^\pm_{(a,b)}$
is constant except when $S$ crosses $\Sigma^{(1)}_{tri}$ in strata of type $(a,b)^\pm$ or $(b,a)^\pm$.
At each such intersection in positive (respectively negative) direction,
$I^\pm_{(a,b)}$ changes exactly by $-1$ (respectively $+1$).
Indeed, the configurations of the three crossings which come together
in a triple point are shown in Fig.~\ref{perturb}.

In each of the four cases, exactly one pair $p, q$ of crossings contributes
to one of the sums $I$.

\begin{figure}
\centering
\includegraphics{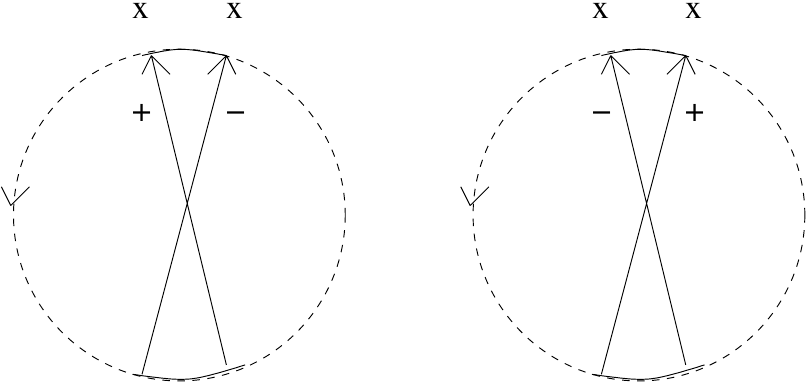}
\caption{\label{autotang} Contributions of auto-tangencies}  
\end{figure}

After drawing all possible triple points, it is easily seen that $p, q$
must verify: $w(p)w(q) = -1$ for the first two cases and
$w(p)w(q)=+1$ for the last two cases.
Notice that the type of the triple point is completely determined by
the sub-configurations shown in Fig.~\ref{subdia}.
Thus, the sums $I$ are constant by passing all types of triple
points except those shown in Fig.~\ref{perturb}. The generic loop $S$ intersects
$\Sigma$ only in strata that correspond to triple points or to
auto-tangencies. An auto-tangency adds to the Gauss diagram always one
of the sub-diagrams shown in Fig.~\ref{autotang}.

The two arrows evidently do not enter together in the configurations shown in Fig.~\ref{subdia}. If one of them contributes to such a configuration, then the other contributes to the same configuration but with an opposite sign.

Therefore, for any two $\hat \beta_1$, $\hat \beta_2 \in M(\hat\beta)$, the difference $I^\pm_{(a,b)}(\hat \beta_1) -
I^\pm_{(a,b)}(\hat \beta_2)$ is just the algebraic intersection
number of an oriented arc from $\hat \beta_1$ to $\hat \beta_2$
with the union of the cycles of codimension one
$(a,b)^+ \cup (b,a)^+$ (resp., $(a,b)^- \cup (b,a)^-$).
Hence, for each loop $S$, these numbers are 0. The identities (*) follow.$\Box$

\begin{example}
Let $\hat \beta$ be the closure of the 4-braid $\beta = \sigma_1
\sigma_2^{-1}\sigma_3^{-1}$. We consider

$\Gamma_{(1,2)^-}$,   $\Gamma_{(2,1)^-}$,
$\Gamma_{(2,3)^+}$ and $\Gamma_{(3,2)^+}$.

A calculation by hand gives:
\begin{displaymath}
\Gamma_{(1,2)^-}(rot(\hat \beta)) =
\Gamma_{(2,3)^+}(rot(\hat \beta)) = -1
\end{displaymath}
and
\begin{displaymath}
\Gamma_{(2,1)^-}(rot(\hat \beta)) =
\Gamma_{(3,2)^+}(rot(\hat \beta)) = +1
\end{displaymath}
Therefore, all four 1-cocycles of Gauss degree 0 are non-trivial and they are obviously related to the finite type invariants  $W_a(\hat \beta)$ of degree 1.
\end{example}

 The non-triviality of a one-cocycle invariant on the loop $rot(\hat\beta)$ implies  that the braid $\beta$ is not periodic, because of the following result of Hatcher.

Let $hat(\hat \beta)$ be the loop which is obtained by the rotation of the solid torus around the core of the complementary solid torus $S^3 \setminus V$. Evidently, $hat(\hat \beta)$ is just a rotation of diagrams in the annulus and hence it does not intersect the discriminant $\sum^{(i)}, i>0,$ at all. Consequently, each one-cocycle invariant vanishes on this loop.

\begin{proposition} (Hatcher).
The loops $rot(\hat \beta)$ and $hat(\hat \beta)$ represent linearly dependent homology classes in $H_1(M(\hat \beta);\mathbb{Q})$ if and only if $ \hat \beta$ is isotopic to a torus knot in $\partial V$. 
\end{proposition}

{\em Proof.}  The rotations give an action of $ \mathbb{Z} \times \mathbb{Z}$, and one can look at the induced action on the fundamental group of the knot complement, with respect to a base point in the boundary torus $\partial V$ of the solid torus.  Call this fundamental group  G.  The action of $ \mathbb{Z} \times \mathbb{Z}$ on G is conjugation by elements of the $ \mathbb{Z} \times \mathbb{Z}$ subgroup of G represented by loops in the boundary torus.  The action will be faithful (zero kernel) if the center of G is trivial.  The only (irreducible) 3-manifolds whose fundamental groups have a nontrivial center are Seifert fibered manifolds defined by a circle action.  The only such manifolds that embed in $\mathbb{R}^3$ and have two boundary components (as here) are the obvious ones that correspond to cable knots, that is, the complement of a knot in a solid torus that is isotopic to a nontrivial loop in the boundary torus.
$\Box$

On the other hand, it is well known that a braid $\beta$, which closes to a knot, is {\em periodic} if and only if $\hat \beta$ is isotopic to a torus knot in $\partial V$ (see e.g. \cite{BGG}).
Consequently, if some one-cocycle evaluated on $rot(\hat \beta)$ is non-trivial, then the braid $\beta$ is not periodic. The above braids are evidently not reducible and as well known this implies now that they are pseudo-Anosov, compare \cite{Th}, \cite{B-B04}. 

It would be very interesting to find out which information about the entropy of a braid $\beta$ and about the simplicial volume of its mapping torus $S^3 \setminus (\hat \beta \cup L)$ can be obtained from the set of all one-cocycle polynomials $\Gamma(rot(\hat \beta))$ (compare the next section). But unfortunately a computer program is still missing in order to calculate lots of examples and to make precise conjectures.

\subsection{One-cocycles of higher Gauss degree}

 A general procedure in order to define all one-cocycle invariants with our method was given in the preprint \cite{F06}. However, it seems to be not very understandable. Therefore we restrict our-self in this paper to the construction of very simple examples, which illustrate our method and which can be generalized easily by an interested reader.

Let $\beta \in B_n$ such that $\hat \beta \hookrightarrow V$ is a knot. We assume that $n$ is divisible by 3.
 We consider only the type of the triple crossings $(2/3n,2/3n)^+$ and all markings in the configurations have to be of type $1/3n$ or $2/3n$.

\begin{figure}
\centering
\includegraphics{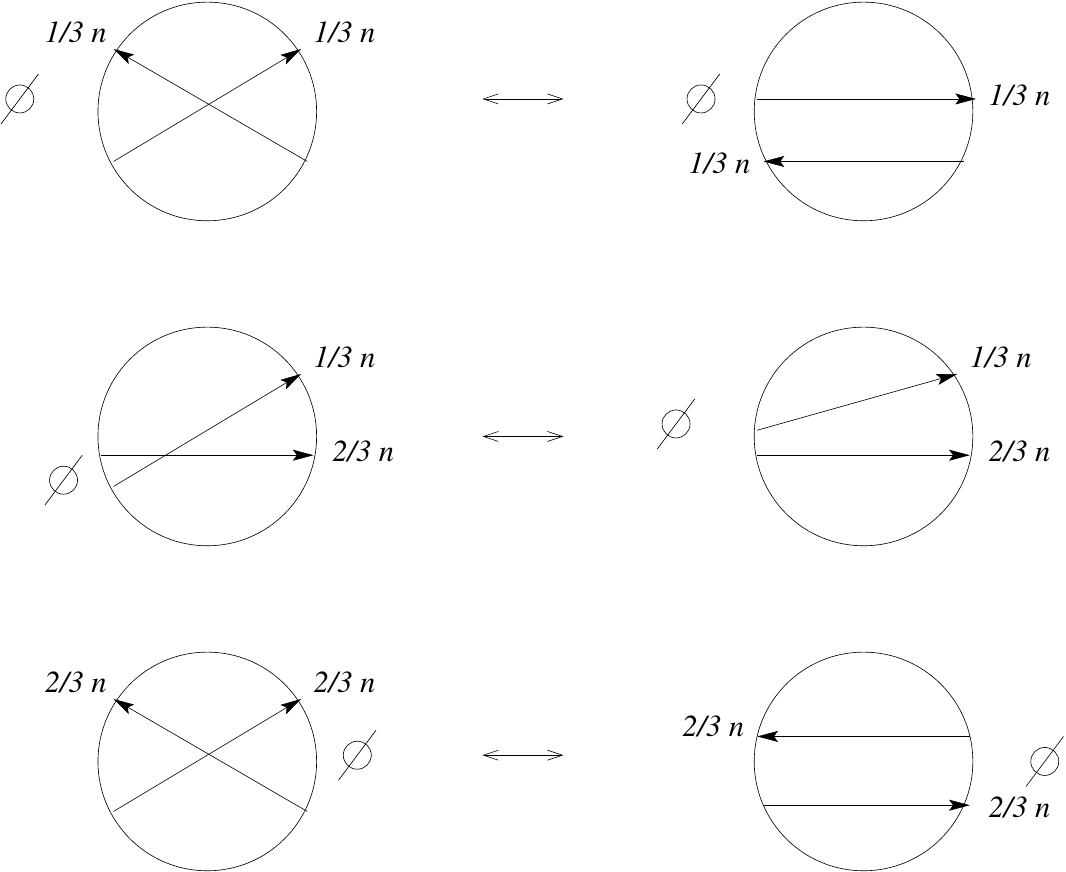}
\caption{\label{slidecup} Slide moves for a couple of arrows not in the triangle}  
\end{figure}

We define a {\em slide move for a couple of arrows} in a configuration in Fig.~\ref{slidecup}. Here $\emptyset$ means that there are no heads or foots of arrows of the configuration in the corresponding segment on the circle.

\begin{figure}
\centering
\includegraphics{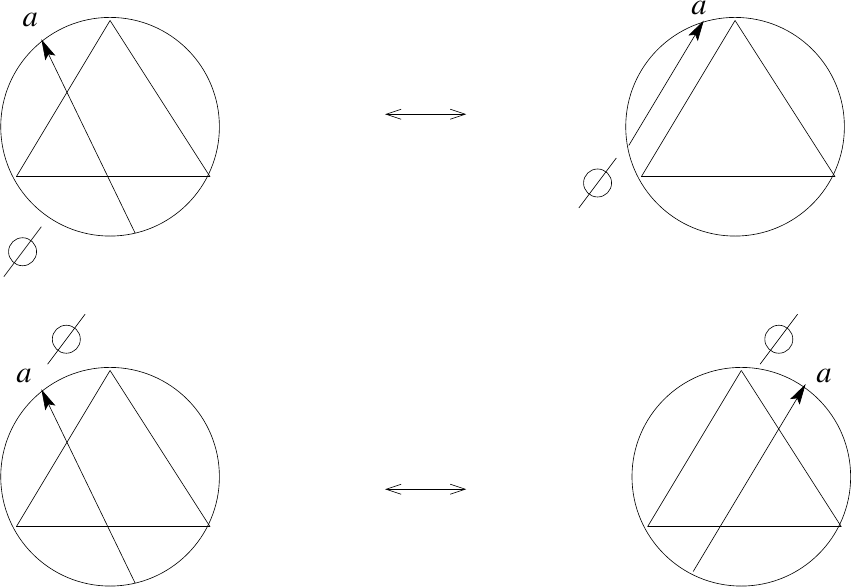}
\caption{\label{slidetri} Slide moves of an arrow with respect to the triangle}  
\end{figure}

We define  {\em slide moves for an arrow with respect to the triangle} in Fig.~\ref{slidetri}.
Here the condition is that the move does not create an oriented loop in the diagram which represents a homology class which is not in $\{ 1 ,2 , \dots , n-1\}$, compare Observation 3. (The marking $a$ here is $1/3n$ or respectively $2/3n$.)

\begin{figure}
\centering
\includegraphics{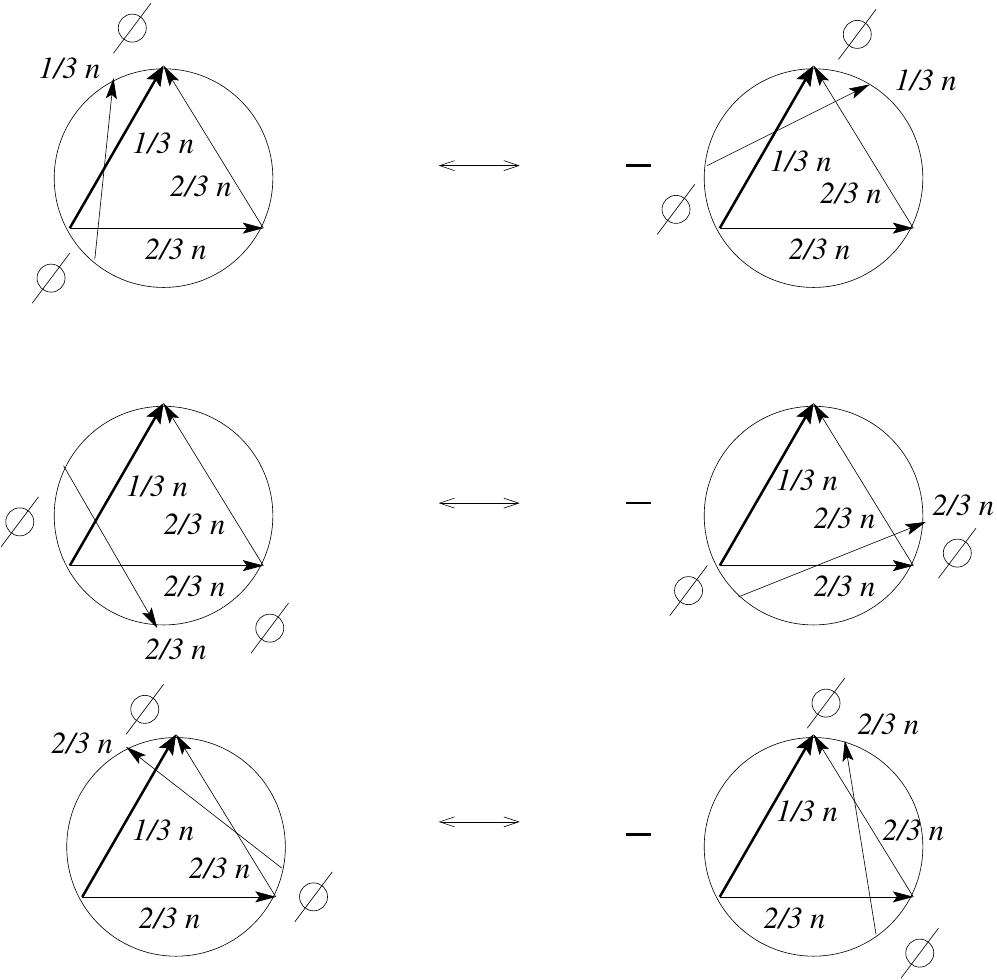}
\caption{\label{exchange} Exchange move with an arrow in the triangle}  
\end{figure}

We define an {\em exchange move for an arrow with the triangle} in Fig.~\ref{exchange}. Notice that the sign of the configuration changes in this case.

\begin{figure}
\centering
\includegraphics{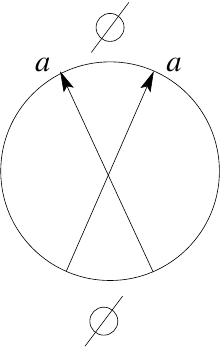}
\caption{\label{forbidden} A forbidden sub-configuration}  
\end{figure}

We define a {\em forbidden sub-configuration} in Fig.~\ref{forbidden}, where the arrows have the same homological marking $a$. (The marking $a$ here is $1/3n$ or respectively $2/3n$.)

\begin{theorem}
Let $I^d= \sum_{i}{\epsilon_i I_i}$ be a linear combination of  configurations of Gauss degree $d$ and which is invariant under all possible slide moves, exchange moves and which does not contain any forbidden sub-configuration. Then   
\begin{displaymath}
\Gamma^d(S) =
\sum_{s_i \in S \cap \Sigma^{(1)}_{tri} \textrm{of type $(2/3n,2/3n)^+$}}
 {sign(s_i)x^{\sum_{i}{\epsilon_i (\sum_{D_i}{\prod_{j}{w(p_j)})}}}}
\end{displaymath}
is a one-cocycle, where $D_i$ is the set of unordered $d$-tuples $(p_1, \dots, p_d)$
of arrows which enter in $I_i$ in the Gauss diagram of $s_i$ .
\end{theorem}

{\em Proof.} Following Theorem 2 and Lemma 3, we have to show that $\Gamma^d$ vanishes on all meridians of $\sum^{(2)} = \sum ^{(2)}_{quad} \cup \sum ^{(2)}_{trans-self} \cup \sum ^{(2)}_{self-flex} \cup \sum^{(2)}_{inter}$ in $M_n$.

The Gauss diagrams in the meridian of $\sum^{(2)}_{inter}$ for two triple crossings change just by slide moves for a couple of arrows and the two triple crossings enter with different signs, compare \cite{F-K}. Consequently  $\Gamma^d=0$ on the meridian. In the transverse intersection of $\sum^{(1)}_{tri}$ with $\sum^{(1)}_{tan}$ there appear two new arrows in the Gauss diagram of one of the two triple crossings in the meridian. These two arrows do not enter together into any configuration because $I^d$ does not contain any forbidden sub-configurations. The two arrows have different signs and consequently, if one of the arrows enter into a configuration then the other enters too but with a different sign of the weight and their contributions cancel out.

Auto-tangencies do not contribute to $\Gamma^d$ and hence $\Gamma^d=0$ on the meridian of $\sum ^{(2)}_{self-flex}$

In the meridian of $\sum ^{(2)}_{quad}$ there are exactly eight triple crossings. They come in pairs with different signs, compare \cite{F-K}. The Gauss diagrams of the pairs differ just by the slide moves of three different arrows with respect to the triangle and hence again
$\Gamma^d=0$ on the meridian, i.e. $\Gamma^d$ is a solution of the tetrahedron equation.

In the meridian of $\sum ^{(2)}_{trans-self}$ there are exactly two triple crossings and they have different signs. The two Gauss diagrams of the braids with the triple crossing differ exactly by an exchange move of one arrow with the triangle. Moreover, the sign of the arrow changes (but we do not care about the signs of the arrows in the triangle), compare \cite{F-K}. But  $I^d$ stays invariant,  because the sign of the configuration changes too for the exchange move with a triangle. Consequently, $\Gamma^d=0$ on the meridian. $\Box$

\begin{figure}
\centering
\includegraphics{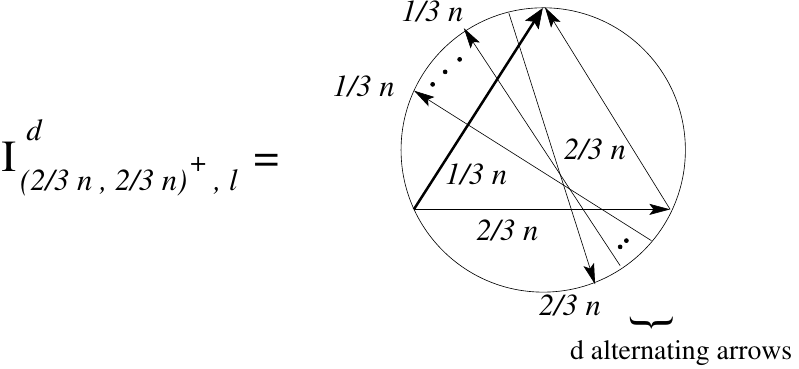}
\caption{\label{Gammal} Configuration for even $d$}  
\end{figure}

\begin{figure}
\centering
\includegraphics{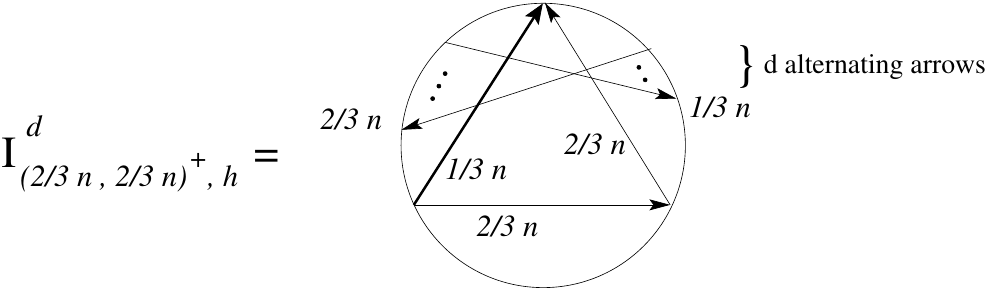}
\caption{\label{Gammah} Another configuration for even $d$}  
\end{figure}

\begin{figure}
\centering
\includegraphics{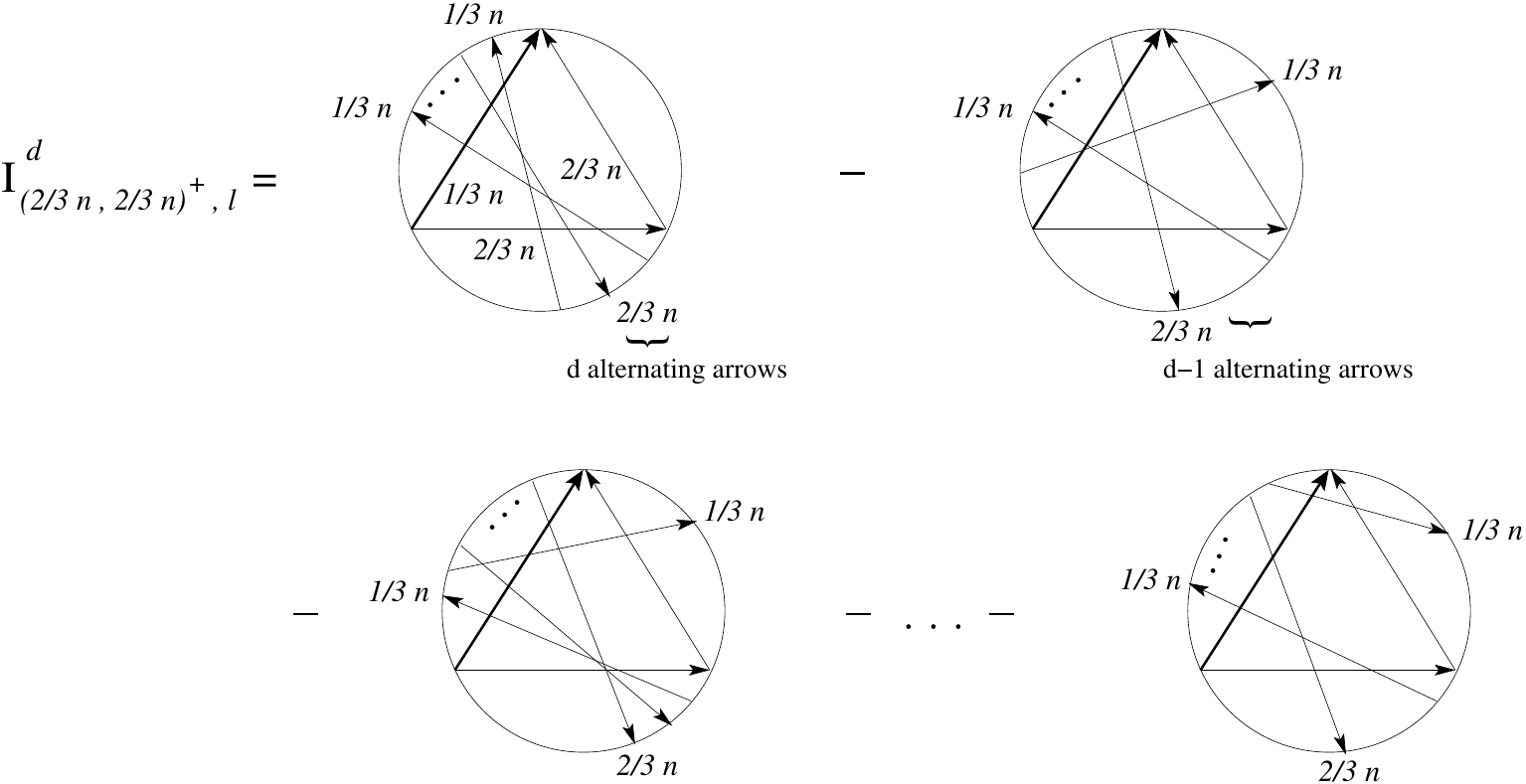}
\caption{\label{Gammaodd} Configuration for odd $d$}  
\end{figure}

We give three examples in Fig.~\ref{Gammal}, Fig.~\ref{Gammah} and Fig.~\ref{Gammaodd}. 

Let us give names to the three strands in a triple crossing: $h$ is the highest strand, $m$ is the strand in the middle and $l$ is the lowest strand. In each of the examples, the bunch of alternating arrows separates one of the strands from the other two. We use the single separated strand in the notation of the configuration (compare the figures).

The three configurations in the examples define one-cocycles because of the following proposition.

\begin{proposition}
The one-cochains $I^d_h$ for even $d$ and $I^d_l$ for all $d$ are invariant under all possible slide moves, exchange moves and do not contain any forbidden sub-configurations. 

\end{proposition}
{\em Proof.} 
$I^d_h$ and $I^d_l$ for even $d$ are {\em rigid}, i.e. no slide move or exchange move at all is possible. Indeed, any slide move would create a loop in the Gauss diagram with homological marking 0. But this contradicts Observation 3. There is no arrow at all which would allow an exchange move with the triangle.
For the same reason, $I^d_l$ for odd $d$ allows only slide moves of just one arrow and one easily sees that $I^d_l$ takes into account all possible slide moves of this arrow, called {\em wandering} arrow, as well as the single possible exchange move for the wandering arrow. $\Box$

But remember that the alternating arrows form only a configuration, i.e. in the Gauss diagram of the closed braid there could be other arrows which cut the alternating arrows in an arbitrary way.

\begin{remark}
The one-cochain which is obtained by taking mirror images of everything in $I^d$, i.e. changing the orientation of each arrow and interchanging the markings $1/3n$ with $2/3n$, is of course also a one-cocycle.

If $n$ is not divisible by 3 then we replace  $\hat\beta$ by a $3k$-cable, $k \in \mathbb{N}^*$, which is twisted by the permutation braid $\sigma_1\sigma_2\ ..\sigma_{3k-1}$ in order to get a knot, and we can take now the markings $kn$ and $2kn$. Indeed, we can imagine the untwisted 3k-cable as $3k$ parallel strands on a band (which projects into the annulus by an immersion)  and we can push the permutation braid along the band. Hence, if two closed braids are isotopic in the solid torus then their 3k-cables which are twisted by the same braid, are still isotopic.
\end{remark}

\begin{example}
Notice that for closed 3-braids, the homological markings of the triple point are already
determined by the arrows. Let us consider the 3-braid $\beta = \sigma_1 \sigma_2^{-1}\sigma_1 \sigma_2 \sigma_1 \sigma_1 \sigma_2 \sigma_1$.
The braid $\beta$ contains a full-twist at the right and we represent $rot$ by pushing it through $\sigma_1\sigma_2^{-1}$ to the left and bringing it back to the right by an isotopy of diagrams in the solid torus. For the convenience of the reader we will give the loop below. As usual we will write shortly $i$ for $\sigma_i$ and $\bar i$ for $\sigma_i^{-1}$ and we put the Reidemeister moves into brackets. 

$1\bar 2 (121)121 \rightarrow 1(\bar 2 2)12121 \rightarrow 1(2 \bar 2)12121 \rightarrow 12(\bar 212)121 \rightarrow 1212(\bar 11)21 \rightarrow 1212( 1\bar 1)21 \rightarrow 12121(\bar 121) \rightarrow 121(212)1 \bar2 \rightarrow 1211211 \bar 2 \rightarrow 1 \bar 2 121121$

Here, the last arrow is an isotopy of diagrams. It turns out that only the second and the forth Reidemeister III move is of type $(2,2)^+,l$. We show the corresponding closed braids together with their Gauss diagrams in Fig.~\ref{Exam}. 

\begin{figure}
\centering
\includegraphics{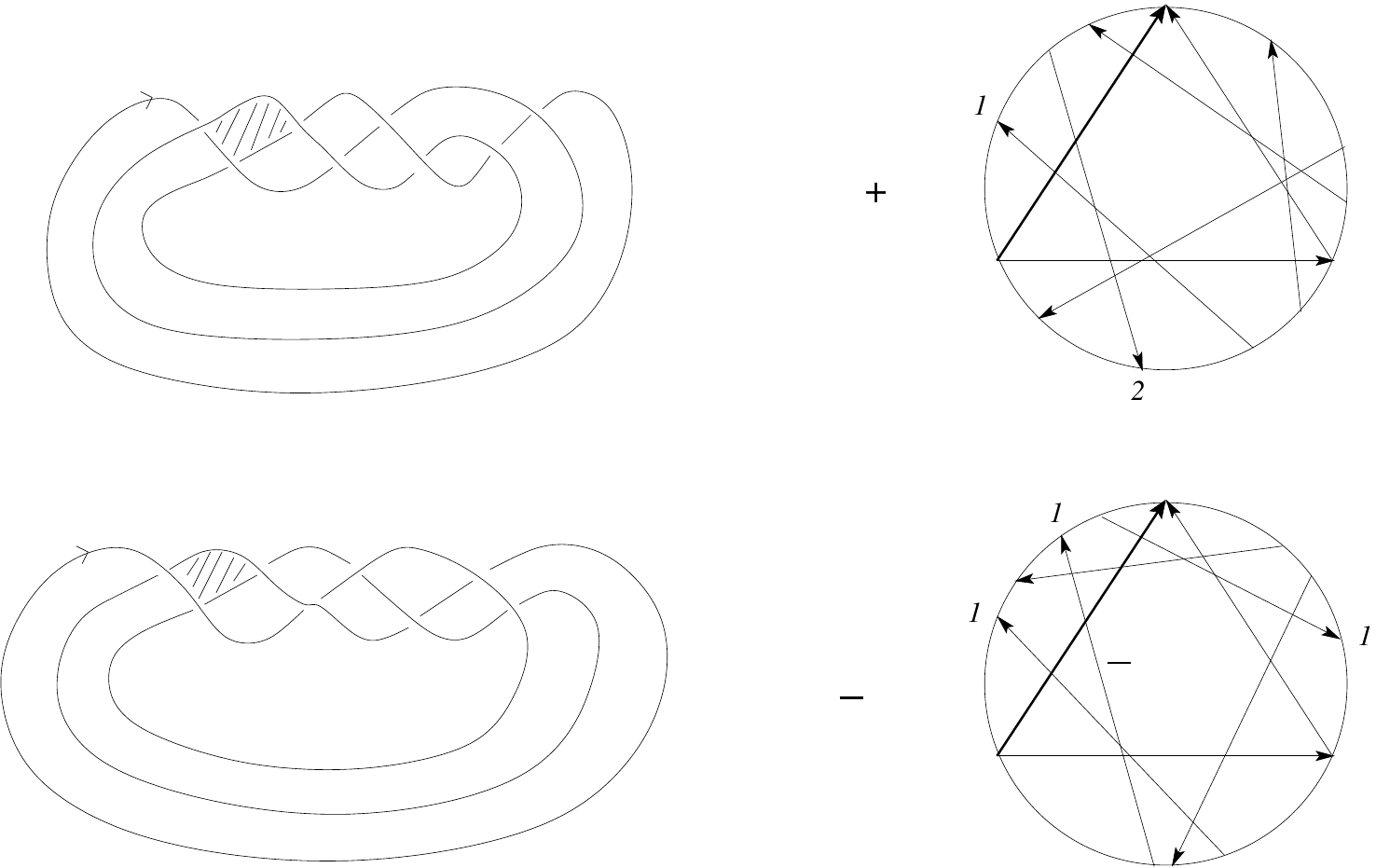}
\caption{\label{Exam} Contributing triple crossings with their Gauss diagrams}  
\end{figure}

We indicate the sign of the move, the homological markings of the relevant crossings and the sign of the crossings, but only  if it is negative. One easily calculates now that

$\Gamma^2_{(2,2)^+,l}(rot(\hat \beta))= x$ and $\Gamma^1_{(2,2)^+,l}(rot(\hat \beta))= x-x^{-1}$

\end{example}

\begin{example}
 Let us consider the closure of the 3-braid $\beta_{+++}= 11122121$. There are exactly eight Reidemeister III moves in pushing $121$ twice through the rest of the closed  braid. A calculation gives now

$\Gamma^4_{(2,2)^+,h}(rot(\hat \beta_{+++}))= -x$, $\Gamma^2_{(2,2)^+,h}(rot(\hat \beta_{+++}))= -x^3-x$, 

$\Gamma^4_{(2,2)^+,l}(rot(\hat \beta_{+++}))= x$, $\Gamma^2_{(2,2)^+,l}(rot(\hat \beta_{+++}))= x^3 + x$, 
$\Gamma^1_{(2,2)^+,l}(rot(\hat \beta_{+++}))= x^2 + x -x^{-1} - x^{-2}$.

Let us consider $\beta_{-++}= 11 \bar 122121$, $\beta_{+-+}= 111\bar 22121$, $\beta_{--+}= 11\bar 1\bar 2 2121$, 
$\beta_{++-}= 1112212\bar1$, $\beta_{-+-}= 11 \bar 12212\bar 1$, $\beta_{+--}= 11 1\bar 2212\bar 1$, $\beta_{---}= 11 \bar 1\bar 2212\bar 1$.

If an invariant $\Gamma (rot(\hat \beta))$ would be of finite type of degree 2 then we would have (compare \cite{BN}, \cite{V90})\vspace{0,2cm}

(*) \hspace{1cm}  $\Gamma (rot(\hat \beta_{+++}))+\Gamma (rot(\hat \beta_{--+}))+\Gamma (rot(\hat \beta_{-+-}))
+\Gamma (rot(\hat \beta_{+--}))-\Gamma (rot(\hat \beta_{-++}))-\Gamma (rot(\hat \beta_{+-+}))-\Gamma (rot(\hat \beta_{++-}))-\Gamma (rot(\hat \beta_{---}))=0$.\vspace{0,2cm}

One easily sees that $\hat \beta_{-++}=\hat \beta_{+-+}=\hat \beta_{++-}=\hat {222121}$ and that $\hat \beta_{--+}=\hat \beta_{+--}=\hat \beta_{-+-}=\hat{1212}$ is the closure of a periodic braid, as well as $\hat \beta_{---}=\hat {12}$. Hence, $(*)$ could be only satisfied if the integer Laurent polynomial $\Gamma (rot(\hat \beta_{+++}))$ is divisible by 3, which is not the case for each of the above polynomials. One can generalize this example to show that the above one-cocycle polynomials are not of finite type of any degree.

\end{example}

\begin{figure}
\centering
\includegraphics{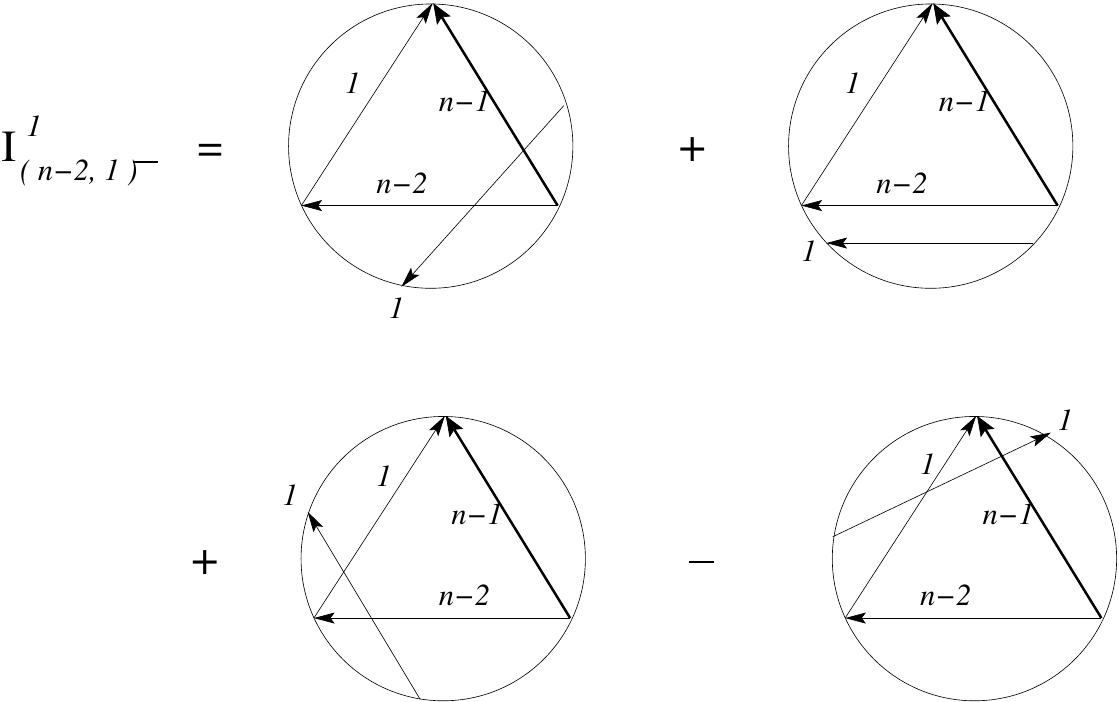}
\caption{\label{Gammafour} A configuration of Gauss degree 1 for $n>3$}  
\end{figure}

Let $n>3$ be arbitrary. We define a linear combination of configurations $I^1_{(n-2,1)^-}=\sum_{i}{\epsilon_i}I^1_i$ in Fig.~\ref{Gammafour}.

\begin{theorem}

\begin{displaymath}
\Gamma^1_{(n-2,1)^-}(S) =
\sum_{s_i \in S \cap \Sigma^{(1)}_{tri} \textrm{of type $(n-2,1)^-$}}
 {sign(s_i)x^{\sum_{i}{\epsilon_i (\sum_{p_i}{w(p_i)})}}}
\end{displaymath}

is a one-cocycle of Gauss degree 1, where $p_i$ are the arrows which enter in $I^1_i$ in the Gauss diagram of $s_i$ .

\end{theorem}

The proof of the theorem is completely analogous to the proof of Theorem 4 and Proposition 4 and is left to the reader.

Again, its mirror image is also a one-cocycle. Moreover, $\Gamma^1_{(n-2,1)^-}$ can be easily generalized to one-cocycles of higher Gauss degrees by introducing appropriate slide moves, along the same lines as Theorem 4 and Proposition 4.

\begin{example}
A calculation yields that for $rot(\hat {\sigma_3\sigma_2\sigma_1^{-1}})$ there is a single triple crossing of type $(2,1)^-$  in the loop $\hat{\beta \Delta^2 \Delta^{-2}} \rightarrow \hat{\Delta^2\beta\Delta^{-2}}$. It is represented by the closure of the 4-braid $(\sigma_2 \sigma_3 \sigma_2) \sigma_2^{-1} \sigma_3^{-2} \sigma_1$, where the parentheses as usual refers to the triple crossing before the Reidemeister III move. One easily calculates now that  
\begin{displaymath}
\Gamma^1_{(2,1)_1 ^-}(rot(\hat {\sigma_1\sigma_2^{-1}\sigma_3^{-1}}))=-x
\end{displaymath}
Therefore, $\Gamma^1_{(2,1)^-}$ is a non-trivial cohomology class for closed pseudo-Anosov  4-braids.
\end{example}

\begin{example}
Let us consider the reducible 4-braid $\beta=\sigma_3\sigma_2\sigma_3\sigma_1\sigma_2$ which is the 2-cable $\sigma_1$ of the 2-braid $\sigma_1$. Evidently, the entropy and the simplicial volume of $\beta$ are trivial, because the JSJ-decomposition of $S^3 \setminus (\hat \beta \cup L)$ has only Seifert fibered pieces.

An easy calculation gives $\Gamma^1_{(n-2,1)^-}(rot(\beta))=x$. Consequently, the vanishing of the entropy and of the simplicial volume of the braid does not imply that the invariant vanishes too.

\end{example}

The 1-cocycle invariants have the following nice property, which can be
used to estimate from below the length of conjugacy classes of braids.
\vspace{0,5cm}

\begin{proposition}
 Let the knot $K = \hat \beta \hookrightarrow V$ be a closed
$n$-braid and let $c(K)$ be its minimal crossing number, i.e. its
minimal word length in $B_n$. Then all 1-cocycle polynomials
of Gauss degree $d$ vanish for
\begin{displaymath}
d \geq c(K) + n^2 - n - 2.
\end{displaymath}
\end{proposition}

{\em Proof:\/}
Assume that the word length of $\beta$ is equal to $c(K)$.
We can represent $[rot(K)]$ by the following
isotopy which uses shorter braids.
\begin{displaymath}
\beta \to \Delta\Delta^{-1}\beta \to \Delta^{-1}\beta\Delta
\to \Delta^{-1}\Delta\beta' \to \beta' \to \Delta\Delta^{-1}\beta'
\to \Delta^{-1}\beta'\Delta \to \Delta^{-1}\Delta\beta \to \beta
\end{displaymath}

Here, $\beta'$ is the result of rotating $\beta$ by $\pi$ i.e.
each $\sigma_i^{\pm 1}$ is replaced by $\sigma_{n-i}^{\pm 1}$.
Obviously, $c(\Delta) = \frac{n(n-1)}{2}$. Thus, each Gauss diagram which
appears in the isotopy has no more than $c(K) + n^2 - n$ arrows. Indeed,
we create a couple of crossings by pushing $\Delta$ through $\beta$
only after having eliminated a couple of crossings before (compare Section 2.2). Therefore, for each diagram with a triple crossing
there are at most $c(K)+n^2-n-3$ other arrows, and hence, each
summand in a 1-cocycle of Gauss degree $d$ is already zero if $d \geq
c(K) + n^2 - n - 2$.$\Box$

(Remember that the degree in Vassiliev's sense of $d(\Gamma(\hat \beta))/dx$ evaluated at  $x=1$ is at most $d+1$, compare Remark 4, which leads to the result cited in the Introduction.)

\section{Character invariants}
In this section we refine our invariants to  another class of easily calculable isotopy invariants for closed braids.

Let $\hat \beta_s , s \in [0 ,1]$ be a generic isotopy of closed braids, $l$ an integer, $l.rot(\hat\beta_s)$ $l$-times the canonical loop and $TL\tilde(l.rot(\hat \beta_s))$ the union of the trace circles of the corresponding resolution of the trace graphs (compare Section 2).
It follows from Observation 1 (compare Remark 2) that the trace circles for different parameter s are in a natural one-to-one correspondence. Consequently, we can give {\em names \/}
$x_i$ to the circles of $TL\tilde(l.rot(\hat \beta_0))$ and extend these names in a unique way on the whole family of trace circles.

Let $\{x_1 ,x_2 ,\dots \}$ be the set of named trace circles. Obviously, for each circle $x_i$ there is a well defined homological marking $h_i \in H_1(V)$.
Let $[x_i] \in H_1(T^2)$ be the homology class represented by the trace circle $x_i$ itself (with its natural orientation induced from the orientation of the trace graph).

\subsection{Character invariants of Gauss degree 0}

In this section, we use the named cycles, i.e. the trace circles, in order to refine the
1-cocycles of Gauss degree 0 which were defined in Section 3.3.

Let $X = \{ x_1, \dots, x_m \}$
be the set of named trace circles. Let $x_{i_1}, x_{i_2}, x_{i_3} \in X$ be fixed.
We do not assume that they are necessarily different. Let $h_{i_1},
h_{i_2}, h_{i_3}$ be the corresponding homological markings.

\begin{figure}
\centering
\includegraphics{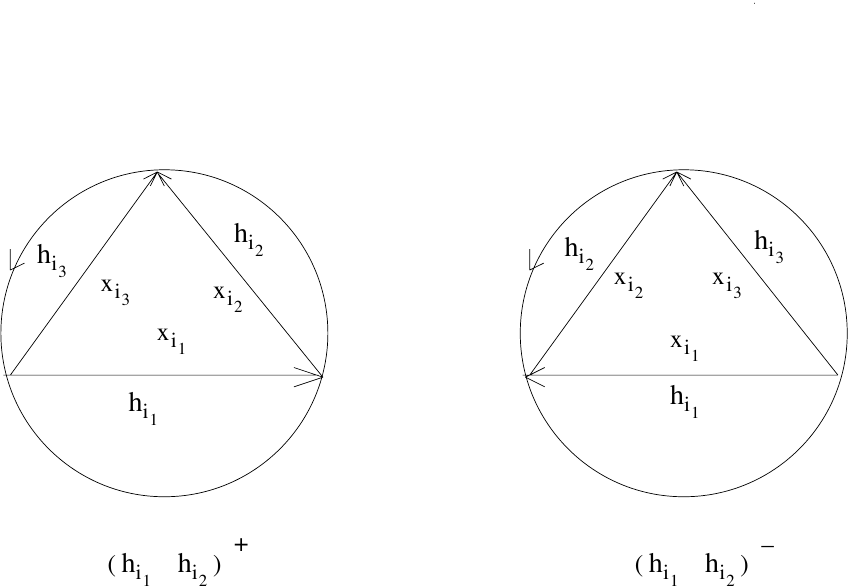}
\caption{\label{nametri} Triangles with names of crossings}  
\end{figure}

\begin{definition}
A {\em character of Gauss degree 0 of $l.rot(\hat\beta)$\/}, denoted by
\begin{displaymath}
C_{(h_{i_1}, h_{i_2})^\pm(x_{i_1}, x_{i_2}, x_{i_3})}(l.rot(\hat\beta))
\end{displaymath}
or sometimes shortly $C(\hat\beta)$, is the algebraic intersection number of
$l.rot(\hat\beta)$ with the strata $(h_{i_1}, h_{i_2})\pm$ in
$\Sigma^{(1)}_{tri}$  and such that the crossings of the triple
point belong to the named trace circles as shown in Fig.~\ref{nametri}.
We call the unordered set $\{ x_{i_1}, x_{i_2}, x_{i_3} \}$ the
{\em support\/} of the character $C(\hat\beta)$.
\end{definition}

\begin{remark}
Evidently, in order to obtain a non-trivial intersection number we
need that $h_{i_3} = h_{i_1} + h_{i_2} - n$ in $(h_{i_1}, h_{i_2})^+$
and that respectively $h_{i_3} = h_{i_1} + h_{i_2}$ in $(h_{i_1}, h_{i_2})^-$
(compare Fig.~\ref{marktype}). 
\end{remark}
Notice that for characters of Gauss degree 0 the relations (*) from
Proposition 2 are no longer valid. For example,
$C_{(h_{i_1}, h_{i_2})^+(x_{i_1}, x_{i_2}, x_{i_3})}(l.rot(\hat\beta))$
can be non-trivial for $h_{i_1} = h_{i_2}$ and even for
$x_{i_1} = x_{i_2}$.

\begin{theorem}
Let $\hat\beta_0$ and $\hat\beta_1$ be isotopic closed braids and let
$\{ x_1, \dots, x_m \}$, $\{ x'_1, \dots, x'_{m'} \}$ be the corresponding
sets of named trace circles of $TL\tilde(l.rot(\hat \beta_0))$ respectively $TL\tilde(l.rot(\hat \beta_1))$. 

Then $m = m'$ and there is a bijection
$\sigma: \{ x_1, \dots, x_m \} \to \{ x'_1, \dots, x'_{m'} \}$ which
preserves the homological markings $h_i$ as well as the homology
classes $[x_i]$ and such that
\begin{displaymath}
C_{(h_{i}, h_{j})^\pm(x_{i_1}, x_{i_2}, x_{i_3})}(l.rot(\hat\beta_0)) =
C_{(h_i, h_j)^\pm(\sigma (x_{i_1}), \sigma (x_{i_2}), \sigma (x_{i_3}))}(l.rot(\hat\beta_1))
\end{displaymath}
for all triples $(x_{i_1}, x_{i_2}, x_{i_3})$.
\end{theorem}

{\em Proof:\/} This is an immediate consequence of Lemma 4
 and the fact that the trace circles are isotopy invariants of $\hat\beta$. $\Box$

\subsection{Character invariants of higher Gauss degree}

We refine the results of the Sections 3.4 and 4.1 in a
straightforward way.

Let $I^d$ be one of the configuration of Gauss degree $d$, which were considered in Section 3.4 and
let $(x_{i_1}, \dots, x_{i_{d+3}})$ be a fixed $(d+3)$-tuple of
elements in $X$ (not necessarily distinct). A {\em named
configuration\/} $I_{(x_{i_1}, \dots, x_{i_{d+3}}), \phi }$ is the
configuration $I^d$ together with a given bijection $\phi$ of
$(x_{i_1}, \dots, x_{i_{d+3}})$ with the $d+3$ arrows in $I^d$ and such
that $x_{i_1}, x_{i_2}, x_{i_3}$ are the arrows of the triangle
exactly as in the previous section. The wandering arrow in $I^d_{(2/3n,2/3n)^+,l}$ has the same name in all configurations and corresponding arrows in the bunches of alternating arrows have of course also the same name. 

We call the corresponding character polynomial $C^d_{(x_{i_1}, \dots, x_{i_{d+3}})}(l.rot(\hat\beta))$.

Different bijections give in general of course different named configurations.

\begin{theorem}
Let $\hat\beta_0$ and $\hat\beta_1$ be isotopic closed braids and let
$\{ x_1, \dots, x_m \}$, $\{ x'_1, \dots, x'_{m'} \}$ be the corresponding
sets of named trace circles of $TL\tilde(l.rot(\hat \beta_0))$ respectively $TL\tilde(l.rot(\hat \beta_1))$. 

Then $m = m'$ and there is a bijection
$\sigma: \{ x_1, \dots, x_m \} \to \{ x'_1, \dots, x'_{m'} \}$ which
preserves the homological markings $h_i$ as well as the homology
classes $[x_i]$ and such that

\begin{displaymath}
C^d_{(x_{i_1}, \dots, x_{i_{d+3}})}(l.rot(\hat\beta_0)) =
C^d_{(\sigma (x_{i_1}), \dots, \sigma(x_{i_{d+3}}))}(l.rot(\hat\beta_1))
\end{displaymath}
for the configurations $I^d$ which were introduced in Section 3.4, and which determine one-cocycles $\Gamma^d$.
\end{theorem}

{\em Proof:\/} This is completely analogous to the proofs of
Theorems 4, 5 and 6. $\Box$

Let $h_i ,h_j$ and a type of triple point, e.g. $(h_i ,h_j)^+$, be fixed. It follows immediately from the definitions that\vspace{0.4cm}

(**)            $\Gamma_{(h_i ,h_j)^+} = \sum_{(x_k ,x_l ,x_m)} C_{(h_i ,h_j)^+}(x_k ,x_l ,x_m )$\vspace{0.4cm}

where $h(x_k) = h_i$ and $h(x_l) = h_j$.

Hence, character invariants define splittings of one-cocycle invariants. However, the set of character invariants on the
right hand side of (**) is not an ordered set. 

It follows from Lemma 2  that the names $x_i$ are determined by their homological markings $h_i$, and that there are exactly $n-1$ trace circles in the case $l=1$.

However, this is in general no longer true in the case of multiples of the canonical loop.

Let $l \in \mathbb{Z}$ be fixed and let $l .rot(\hat \beta)$ be the loop which is defined by going $l$ times along the canonical loop.
Let $TL(l.rot(\hat\beta))$ denote the corresponding trace graph (the t-coordinate in the thickened torus covers now $l$ times the t-circle).

We will show in a simple example that the monodromy is in general non-trivial already for $l =2$.

\begin{example}
Let $\beta = \sigma_2\sigma_1^{-1}\sigma_2\sigma_1^{-1} \in B_3$.
We will write shortly $\beta = 2\bar12\bar1$.

The combinatorial canonical loop for $l = 2$ is given by the following sequence (where we write the names of the crossings just below the crossings).
$2\bar 12\bar 1 \hspace{0.3cm} \to \hspace{0.3cm}
\bar 1\bar 2\bar 12\bar 12(\bar 11)21 \hspace{0.3cm} \to \hspace{0.3cm}
\bar 1\bar 2\bar 12\bar 1221 \hspace{0.3cm} \to \hspace{0.3cm}
\bar 1\bar 2\bar 12\bar 121(\bar 121) \hspace{0.4cm} *_1\to$\\\vspace{0.1cm}
\hspace*{0.5cm} $abcd \hspace{1.3cm}
zyxabcdxyz \hspace{1.4cm}
zyxabcyz \hspace{1.5cm}
zyxabcu_1u_1yz$\\
$\bar 1\bar 2\bar 12\bar 1(212)1\bar 2 \hspace{0.3cm} *_2\to
\bar 1\bar 2\bar 12(\bar 11)211\bar 2 \to
\bar 1\bar 2\bar 12211\bar 2 \to
\bar 1\bar 2\bar 121(\bar 121)1\bar 2 \hspace{0.3cm} *_3\to$
\vspace{0.1cm}
\hspace*{0.5cm} $zyxabcu_1zyu_1 \hspace{1.2cm}
                 zyxabzu_1cyu_1 \hspace{0.5cm}
                 zyxau_1cyu_1 \hspace{0.5cm}
                 zyxau_2u_2u_1cyu_1$\\

$\bar 1\bar 2\bar 1(212)1\bar 21\bar 2 \hspace{0.5cm} *_4\to
(\bar 1\bar 2\bar 1121)1\bar 21\bar 2 \hspace{0.3cm} \to \hspace{0.3cm}
1\bar 21\bar 2 \hspace{0.2cm} \to \hspace{0.2cm}
\bar 1\bar 2\bar 11\bar 21(\bar 212)1 \hspace{0.2cm} *_5\to$
\vspace{0.1cm}
\hspace*{0.2cm} $zyxau_2cu_1u_2yu_1 \hspace{0.9cm}
                 zyxcu_2au_1u_2yu_1 \hspace{0.6cm}
                 u_1u_2yu_1 \hspace{0.5cm}
                 z_1y_1x_1u_1u_2yu_1x_1y_1z_1$\\

$\bar 1\bar 2\bar 11\bar 2112(\bar 11) \hspace{0.3cm} \to \hspace{0.3cm}
\bar 1\bar 2\bar 11\bar 2112 \hspace{0.3cm} \to \hspace{0.3cm}
\bar 1\bar 2\bar 11\bar 21(121)\bar 1 \hspace{0.3cm}*_6\to \hspace{0.3cm}
\bar 1\bar 2\bar 11(\bar 212)12\bar 1$
\vspace{0.1cm}
$z_1y_1x_1u_1u_2yy_1x_1u_1z_1 \hspace{0.2cm}
z_1y_1x_1u_1u_2yy_1x_1 \hspace{0.2cm}
z_1y_1x_1u_1u_2yy_1x_1v_1v_1 \hspace{0.3cm}
z_1y_1x_1u_1u_2yv_1x_1y_1v_1$\\

$*_7\to
\bar 1\bar 2\bar 1112(\bar 11)2\bar 1 \to
\bar 1\bar 2\bar 11122\bar 1 \to
\bar 1\bar 2\bar 11(121)\bar 12\bar 1 \hspace{0.3cm} *_8\to
(\bar 1\bar 2\bar 1121)2\bar 12\bar 1$
\vspace{0.1cm}
$z_1y_1x_1u_1v_1yu_2x_1y_1v_1 \hspace{0.4cm}
                 z_1y_1x_1u_1v_1yy_1v_1 \hspace{0.3cm}
                 z_1y_1x_1u_1v_1yv_2v_2y_1v_1 \hspace{0.3cm}
                 z_1y_1x_1u_1v_2yv_1v_2y_1v_1$\\

$\to 2\bar 12\bar 1$\\
\hspace*{0.5cm} $v_1v_2y_1v_1$

\vspace{0.5cm}

We have the identifications: $d=x$, $b=z$, $x=c$, $y=u_2$, $z=a$,
$u_1=z_1$, $u_2=x_1$, $x_1=u_1$, $y_1=v_2$, $z_1=y$.
This gives us:

$2\bar 12\bar 1 \to 2\bar 12\bar 1$\\
\hspace*{0.5cm} $aacc  \hspace{0.6cm} v_1v_2v_2v_1$

together with the names $u_1=z_1=x_1=u_2=y$.
The second rotation gives us:

$2\bar 12\bar 1 \to 2\bar 12\bar 1 \to 2\bar 12\bar 1$\\
\hspace*{0.5cm} $aacc \hspace{0.4cm} v_1v_2v_2v_1 \hspace{0.4cm}
v_3v_4v_4v_3$

together with the identification $v_1=v_2$, $y_3=v_4$, $z_2=v_2$
and with the names $u_1=z_1=x_1=u_2=y$, $u_3=z_3=x_3=u_4=y_2$.
The monodromy (i.e. how the set of crossings is mapped to itself after the rotation) implies now: $a=v_3=v_4=c$.

Therefore, we have exactly four named cycles: $a, v_1, u_1, u_3$ for

$2\bar 12\bar 1 \to 2\bar 12\bar 1$ with $l = 2$.

Consequently, we have 

$2\bar12\bar1 \to 2\bar12\bar1 \to 2\bar12\bar1$ with the names $aaaa \to v_1v_1v_1v_1 \to aaaa$.

 Hence, $rot(\hat \beta)$ 
acts by interchanging $a$ and $v_1$ (as well as $u_1$ and $u_3$). Consequently, the monodromy is non-trivial in this example.
\end{example}

\subsection{Examples for character invariants of trivial Gauss degree }

In $l.rot(\hat\beta)$ each crossing gives rise to exactly $2l(n-2)$ triple crossings (or Reidemeister III moves). Consequently, the calculation of 

$C_{(h_{i}, h_{j})^\pm(x_{i_1}, x_{i_2}, x_{i_3})}(l.rot(\hat\beta))$ is of linear complexity with respect to the braid length of $\beta$ for fixed $l$ and $n$. Evidently, $n<c+2$ for a braid which closes to a knot, and hence the invariant is of quadratic complexity with respect to the braid length  $c$ for fixed $l$ if $n$ is not fixed.

The following examples for character invariants of Gauss degree 0 are calculated by Alexander Stoimenow using his program in c++. His program is available by request (see \cite{S}).

Let $\beta = \bar12\bar1^32^3 \in B_3$.
($\hat \beta$ represents the knot $8_9$ in the Rolfsen Table.)

We want to show that the link $\hat \beta  \cup$ (core of complementary solid torus) is not invertible in $S^3$.
This is equivalent to show that $\beta$ is not conjugate to $\beta_{inverse} = 2^3\bar1^32\bar1$, i.e. reading the braid backwards (compare e.g. \cite{F01}).

Because $\beta$ is a 3-braid , the homological markings are in $\{1 ,2\}$.

Character invariants of degree one for $l = 1$ and $l = 2$ do not distinguish $\hat \beta$ from $\hat \beta_{inverse}$.
However, for $l = 3$ we obtain three different named cycles  $x_1 ,x_2 ,x_3$ of homological marking 1 and three different named cycles 
$y_1 ,y_2 ,y_3$ of marking 2.

\begin{figure}
\centering
\includegraphics{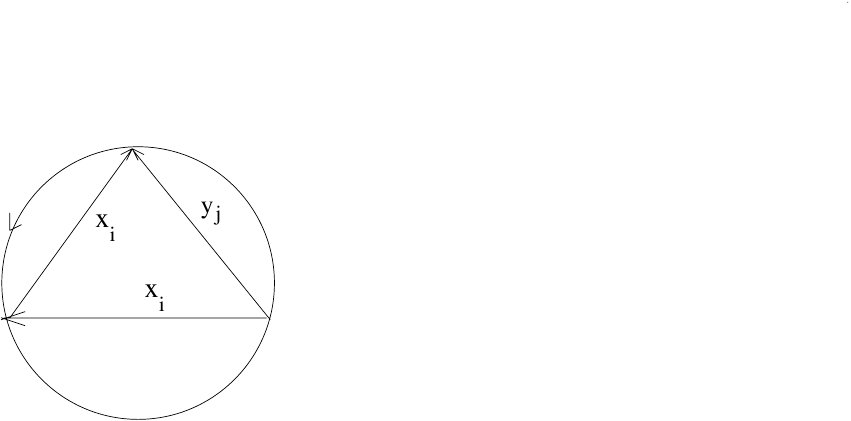}
\caption{\label{carex} The names in the first example}  
\end{figure}

We consider the set of nine character invariants of Gauss degree 0 which are of the form which is shown in Fig.~\ref{carex} , where $i ,j \in \{1 ,2 ,3\}$.

For $\hat \beta$  we obtain the set $\{-1 ,-1 ,-1 ,-1 ,-1 ,-1 ,2 ,2 ,2 \}$ and for $\hat \beta_{inverse}$ we obtain the set 
$\{1 ,1 ,1 ,1 ,1 ,1 ,-2 ,-2 ,-2\}$. Evidently, there is no bijection of the trace circles for $\hat \beta$ and those for 
$\hat \beta_{inverse}$ which identifies the above sets.
Consequently, $\hat \beta$ and $\hat \beta_{inverse}$ are not isotopic in the solid torus. 

The knot $9_5$ can be represented as a 8-braid with 33 crossings. Character invariants of Gauss degree 0 for $l = 2$ show that the braid is not invertible in the same way as in the previous example.

The knot $8_6$ can be represented as a 5-braid with 14 crossings. Character invariants of Gauss degree 0 for $l = 2 ,4 ,6$
show that it is not invertible as a 5-braid. (Surprisingly, it does not work for for $l = 1 ,3 ,5$.)

The knot $8_{17}$ is not invertible as a 3-braid, which is shown with $l = 4$. (It does not work with $l = 1 ,2 ,3$.)

Let $b \in P_5$ be Bigelow's braid (see \cite{Bi}). It has trivial Burau representation. Let $s = \sigma_1\sigma_2\sigma_3\sigma_4$.
The braids $s$ and $bs$ have the same Burau representation. This is still true for their 2-cables, i.e. we replace each strand by two parallel strands.
Character invariants of Gauss degree 0 for $l = 2$ show that the (once positively half-twisted) 2-cables of the above braids are not conjugate, and, consequently, the braids $s$ and $bs$, which have the same Burau representation, are not conjugate either. This shows, that even our character invariants of Gauss degree 0 can not be extracted from the Burau representation.\vspace{0,5cm}

We will construct a further refinement of character invariants of trivial Gauss degree. The number of triple points in a trace graph can only change by trihedron moves, as follows from Theorem 3 .

\begin{definition}
A {\em generalized trihedron \/} is a trihedron which might have other triple points on the edges.
\end{definition}

\begin{figure}
\centering
\includegraphics{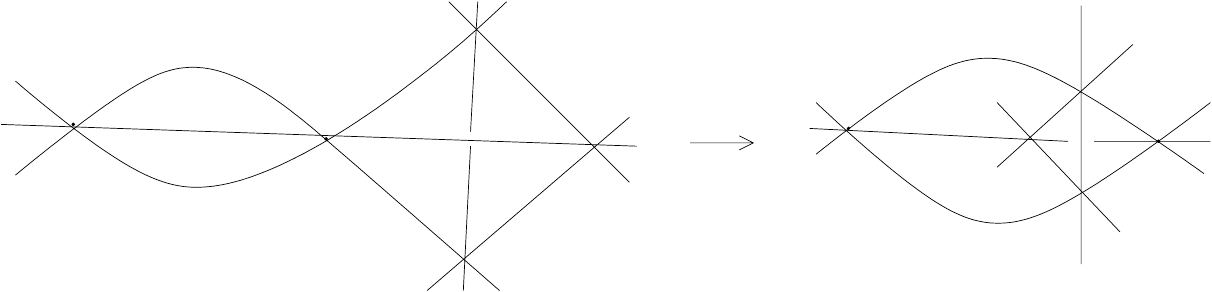}
\caption{\label{gentrihed} A tetrahedron move which creates a generalized trihedron}  
\end{figure}

Fig.~\ref{gentrihed} shows a tetrahedron move which transforms a trihedron into a generalized trihedron. The generalized trihedron has still exactly two vertices.

Evidently, the number of generalized trihedrons does not change under tetrahedron moves. Let $E$ be the set of all triple
points in the trace graph $TL(l.rot(\hat \beta))$ which are {\em not} vertices of generalized trihedrons. The following lemma is an immediate consequence of Theorem 3.

\begin{lemma}
The set $E$, and, hence $card(E)$, is an isotopy invariant of closed braids $\hat \beta$.
\end{lemma}

Moreover, for each element of $E$ we have the additional structure defined before: {\em type , sign, markings, names}.

It follows from Theorem 3 and the geometric interpretation of generalized trihedrons in \cite{F-K} that the two vertices of a generalized trihedron have always different signs. Consequently, character invariants of Gauss degree 0 count just the algebraic number of elements in $E$ which have a given type and given names. But already the geometric number of such elements in $E$ is an invariant as shows Lemma 5.

\begin{definition}
Let $C^{+(-)}_{(h_i ,h_j)^{+(-)}}(x_k ,x_l ,x_m)$ be the number of all positive (respectively negative) triple points in $E$
of given type $(h_i ,h_j)^{+(-)}$ and with given names $x_k ,x_l ,x_m $. We call these the {\em positive (respectively negative) character invariants\/}.
\end{definition}  

The following proposition is now an immediate consequence of Lemma 5 and Definition 15.

\begin{proposition}
The positive and the negative character invariants are isotopy invariants of closed braids for each fixed $l$.
\end{proposition}

\begin{example}
Let us consider $\beta = \sigma_2\sigma_1^{-1} \in B_3$. Its trace graph $TL(rot(\hat \beta))$ is shown in Fig.~\ref{tracegraph}.
One easily sees that it does not contain any generalized trihedrons. Consequently, all four triple points are in $E$.
There are exactly two names $x_1$ and $x_2$. They correspond to the homological markings
$h_1 = 1$ and $h_2 = 2$.

\begin{figure}
\centering
\includegraphics{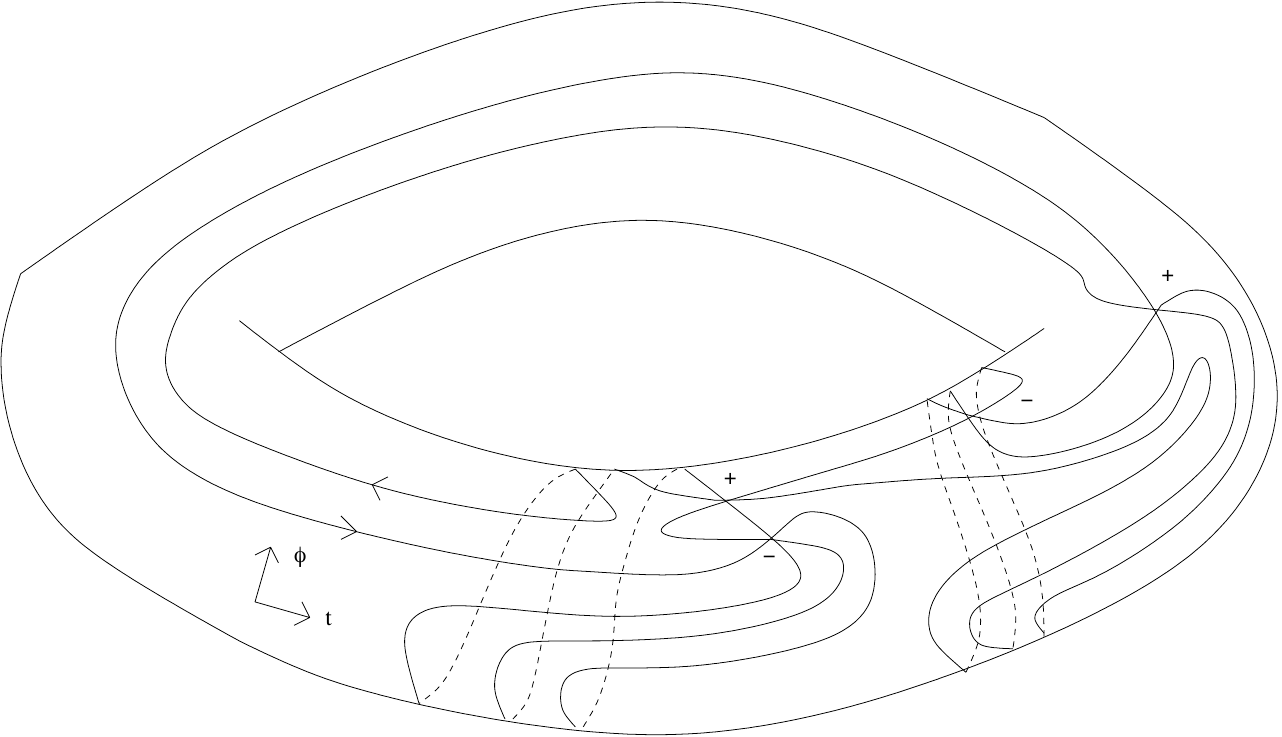}
\caption{\label{tracegraph} The trace graph of $rot(\hat{\sigma_2\sigma^{-1}})$}  
\end{figure}

One easily calculates that two of the triple points are of type $(1 ,1)^-$ and they have different signs. The other two are of type 
$(2 ,2)+$ and they have different signs too. Consequently, all character invariants of Gauss degree 0 are zero.

However, we have $C^+_{(1 ,1)^-}(x_1 ,x_1 ,x_2) = C^-_{(1 ,1)^-}(x_1 ,x_1 ,x_2) = 1$, and\\$C^+_{(2 ,2)^+}(x_2 ,x_2 ,x_1) = C^-_{(2,2)^+}(x_2 ,x_2 ,x_1) = 1$.

Consequently, the positive and negative character invariants contain in this example  for $l=1$ more information than the character invariants of Gauss degree 0.
\end{example}

Unfortunately, there is not yet a computer program available in order to calculate these invariants  as well as character invariants of higher Gauss degree in  more sophisticated examples.

Institute de Math\'ematiques de Toulouse, UMR 5219

Universit\'e Paul Sabatier

118, route de Narbonne 

31062 Toulouse Cedex 09, France

fiedler@math.univ-toulouse.fr
\end{document}